\documentclass[11pt, a4paper,oneside, reqno]{amsart}
\usepackage[english]{babel}
\usepackage{amsmath, amsthm, amsfonts, mathrsfs, amssymb, amscd}
\usepackage{bm}
\usepackage{mathtools}
\usepackage{xcolor}
\definecolor{grey}{rgb}{0.5,0.5,0.5}
\mathtoolsset{centercolon}
\usepackage{accents}
\usepackage{cancel}

\usepackage[toc]{appendix}

\usepackage[shortlabels]{enumitem}
\usepackage[backgroundcolor=yellow, colorinlistoftodos,prependcaption,textsize=small,textwidth=25mm]{todonotes}
\usepackage[colorlinks, citecolor = blue, urlcolor={red}]{hyperref}
\usepackage{geometry}
\allowdisplaybreaks

\newtheorem{theorem}{Theorem}

\newtheorem{lemma}[theorem]{Lemma}

\newtheorem{proposition}[theorem]{Proposition}

\theoremstyle{remark} 
\newtheorem{remark}[theorem]{Remark}

\theoremstyle{definition} 
\newtheorem{definition}[theorem]{Definition}

\numberwithin{theorem}{section}
\numberwithin{equation}{section}

\def\R{{\mathbb R}}

\newcommand{\om}{\omega}

\newcommand{\one}{{{\bf 1}}}

\newcommand{\Dir}{{\rm Dir}}
\newcommand{\Neu}{{\rm Neu}}

\newcommand{\bc}{\mathsf{bc}}
\newcommand{\ii}{{\rm i}} 
\renewcommand{\mod}{{\rm mod\,}}

\newcommand{\mc}{\mathcal}

\newcommand{\Hinf}{H^{\infty}}
\newcommand{\RR}{\mathbb{R}}

\newcommand{\CC}{\mathbb{C}}
\newcommand{\NN}{\mathbb{N}}

\newcommand{\OO}{\mathcal{O}}

\renewcommand{\SS}{\mathcal{S}}

\newcommand{\half}{\frac{1}{2}}
\newcommand{\RRdh}{\RR^d_+}
\newcommand{\RRd}{\RR^d}
\newcommand{\Cc}{C_{\mathrm{c}}}
\renewcommand{\d}{\partial}
\newcommand{\del}{\Delta}

\newcommand{\delDir}{\del_{\operatorname{Dir}}}
\newcommand{\delNeu}{\del_{\operatorname{Neu}}}

\newcommand{\eps}{\varepsilon}

\newcommand{\gam}{\gamma}
\newcommand{\tgam}{\tilde{\gamma}}

\newcommand{\Tr}{\operatorname{Tr}}
\newcommand{\ext}{\operatorname{ext}}

\renewcommand{\tilde}[1]{\widetilde{#1}}
\renewcommand{\hat}[1]{\widehat{#1}}
\newcommand{\dist}{\operatorname{dist}}

\DeclareMathAlphabet{\mathpzc}{OT1}{pzc}{m}{it}

\newcommand{\dd}{\hspace{2pt}\mathrm{d}}

\DeclareMathOperator{\UMD}{UMD}

\author[Floris Roodenburg]{Floris B. Roodenburg}
\address{Floris Roodenburg \hfill\break\indent
Delft Institute of Applied Mathematics \hfill\break\indent
Delft University of Technology \hfill\break\indent
P.O. Box 5031 \hfill\break\indent
2600 GA Delft, The Netherlands}
\email{f.b.roodenburg@tudelft.nl}

\begin{document}
\title[]{Optimal semigroup estimates and functional calculus for the Laplacian on weighted Sobolev spaces}

\makeatletter
\@namedef{subjclassname@2020}{\textup{2020} Mathematics Subject Classification}
\makeatother


\thanks{The author is supported by the VICI grant VI.C.212.027 of the Dutch Research Council (NWO)}

\begin{abstract}
In this paper, we consider the Laplace operator on the half-space with Dirichlet and Neumann boundary conditions. These operators are studied on Sobolev spaces with power weights measuring the distance to the boundary. We prove optimal estimates for the resolvent operators and the corresponding heat semigroups. In addition, it is proved that the Dirichlet and Neumann Laplacians admit a bounded $\Hinf$-functional calculus on Sobolev spaces with certain compatibility conditions at the boundary. We show that these compatibility conditions cannot be omitted in general. The results in this paper are a direct extension of those obtained by Lindemulder, Lorist, the author, and Veraar in [J. Funct. Anal., 289(8):110985, 2025].
\end{abstract}

\maketitle

\setcounter{tocdepth}{1}
\tableofcontents

\section{Introduction} \label{sec:intro}
We establish boundedness of functional calculus for the Dirichlet and Neumann Laplacian and optimal growth bounds of the corresponding Dirichlet and Neumann heat semigroups on power-weighted Sobolev spaces. To be precise, for $p\in (1,\infty)$, $\gam>-1$ and $X$ a Banach space, we consider the weighted Lebesgue space $L^p(\RRdh, w_\gam;X)$ containing all strongly measurable $f:\RRdh\to X$ such that
\begin{equation*}
    \|f\|_{L^p(\RRdh, w_\gam;X)}:= \Big(\int_{\RRdh} \|f(x)\|_X^p w_{\gam} (x)\dd x\Big)^{\frac{1}{p}}<\infty,
\end{equation*}
where 
\begin{equation*}
    w_\gam(x):=\dist(x,\d\RRdh)^\gam=|x_1|^\gam,\qquad x=(x_1, \tilde{x})\in \RR_+\times \RR^{d-1}=\RRdh. 
\end{equation*}
Additionally, for $k\in \NN_0$ define the weighted Sobolev space
\begin{equation*}
    W^{k,p}(\RRdh, w_{\gam};X):=\{f\in \mc{D}'(\RRdh;X): \forall |\alpha|\leq k, \d^\alpha f\in L^p(\RRdh, w_\gam;X)\}
\end{equation*}
equipped with the canonical norm.

The motivation for studying differential operators and their semigroups on these weighted Sobolev spaces stems from the theory of (stochastic) partial differential equations, where the power weights can be used to damp certain irregular behaviour close to the boundary of the domain. In addition, boundedness of the functional calculus plays a prominent role for the study of linear parabolic stochastic partial differential equations. We refer to the introduction of \cite{LLRV24} for an elaborate discussion and overview of the topic. In particular, we note that related results on heat semigroups on weighted spaces, but in different settings, can be found in \cite{Kr99c,MT23}.\\

In the classical unweighted $L^p$-setting, it is well-known that the Dirichlet Laplacian $-\delDir$ and the Neumann Laplacian $-\delNeu$ are sectorial, generate a bounded analytic $C_0$-semigroup, and have a bounded functional calculus. In the case of higher-order Sobolev spaces, one has to additionally impose certain compatibility conditions to obtain a sectorial operator. These unnatural compatibility conditions were shown to be necessary in \cite{DD11}. Nonetheless, without compatibility conditions one can obtain certain $\lambda$-dependent estimates for the resolvent, see, e.g., \cite[Chapter 9]{Kr96_book}.

The Dirichlet and Neumann Laplacian on weighted Sobolev spaces have been studied in \cite{LLRV24, LV18}. Interestingly, it was found by Lindemulder, Lorist, the author, and Veraar in \cite{LLRV24} that the Dirichlet and Neumann Laplacian are \emph{not} sectorial on certain weighted Sobolev spaces and actually generate an analytic $C_0$-semigroup that has polynomial growth. This implies that only the shifted operators $\lambda-\delDir$ and $\lambda-\delNeu$ are sectorial for $\lambda>0$. In particular, almost optimal estimates on the growth of the Dirichlet and Neumann heat semigroup have been established in \cite{LLRV24} for a certain range of exponents $\gam$ in the power weight. 

The aim of the present paper is to sharpen the growth rates for these heat semigroups obtained in \cite{LLRV24} to the optimum. Moreover, we complement the results of \cite{LLRV24} by extending the range of admissible weight exponents. 
In particular, we prove that $\lambda-\delDir$ and $\lambda-\delNeu$ with $\lambda>0$ are sectorial and have a bounded functional calculus on weighted Sobolev spaces with certain compatibility conditions. In addition, we show that these compatibility conditions are necessary to obtain a sectorial operator. We state the main theorem of this paper below.\\

We first introduce the Dirichlet and Neumann semigroups. For $\om\in (0,\pi)$, let $\Sigma_\om=\{z\in \CC\setminus\{0\}:|\arg(z)|< \om\}$ be a sector in the complex plane. Let $G^d_z: \RRd\to \CC$ be the standard heat kernel on $\RRd$ given by 
\begin{equation*}
    G^d_z(x): =\frac{1}{(4\pi z)^{d/2}}e^{-\frac{|x|^2}{4z}},\qquad  z\in \Sigma_{\frac{\pi}{2}}.
\end{equation*}
The standard heat semigroup on $\RR^d$ is given by $T^{d}(z)f(x)=G_z^d\ast f(x)$ for any $f$ such that the convolution is well defined.
The Dirichlet and Neumann heat semigroups on $\RRdh=\RR_+\times \RR^{d-1}$ are obtained by an odd and even reflection, respectively. For $z\in \Sigma_{\frac{\pi}{2}}$ and $x =(x_1, \tilde{x})\in \RRdh$ and $y =(y_1, \tilde{y})\in \RRdh$, we define the kernels
\begin{equation*}
    H^{d,\pm}_z(x,y):=G_z^d(x_1-y_1, \tilde{x}-\tilde{y})\pm G_z^d(x_1+y_1, \tilde{x}-\tilde{y}).
\end{equation*}
Then the Dirichlet and Neumann semigroups $T^d_\Dir$ and $T^d_\Neu$ are given by
\begin{align*}
    T^d_\Dir(z)f(x)& = \int_{\RRdh}H^{d, -}_z(x, y)f(y)\dd y,\\
    T^d_\Neu(z)f(x)&= \int_{\RRdh}H^{d, +}_z(x, y)f(y)\dd y,
\end{align*}
for any $f$ such that the above formulas are well defined. \\

For $p\in (1,\infty)$, $k\in \NN_0$, $\gam>-1$, $\bc\in \{\Dir, \Neu\}$, $k_\Dir:=0$ and $k_\Neu:=1$, we define the weighted Sobolev spaces with compatibility conditions
\[
\begin{aligned}
 W_{\del,\bc}^{k,p}(\RRdh,w_\gamma;X)
 :=\big\{f\in W^{k,p}(\RRdh,w_\gamma;X):\Tr   \d_1^{k_\bc+2j}f=0,\,
 \forall j\in\NN_0:  k_\bc+2j <k-\tfrac{\gamma+1}{p}\big\}.
\end{aligned}
\]
All the traces in this definition are well defined, see \cite[Section 3.2]{LLRV24}. In particular, if $\gam>(k-k_\bc)p-1$, then
\[
    W_{\del, \bc}^{k,p}(\RRdh,w_\gamma;X)
    =W^{k,p}(\RRdh,w_\gamma;X)
\]
and if $(k-k_\bc)p-1<\gam < (k+2-k_\bc)p-1$, then 
\[
    W_{\del, \bc}^{k+2,p}(\RRdh,w_\gamma;X)
    =W^{k+2,p}_{\bc}(\RRdh,w_\gamma;X):= \big\{f\in W^{k+2,p}(\RRdh,w_\gamma;X):\Tr\d_1^{k_\bc} f =0\big\}.
\]

The main result of this paper reads as follows. We refer to Section \ref{sec:elliptic} for the notions of sectoriality and bounded $\Hinf$-calculus. For the definition of $\UMD$ Banach spaces, we refer to \cite[Chapter 4]{HNVW16}. In particular, we note that $\CC$ is a $\UMD$ Banach space.

\begin{theorem}\label{thm:intro_main_thm}
    Let $p\in (1,\infty)$, $k\in \NN_0$, $\bc\in \{\Dir, \Neu\}$, $\gam\in (-1, (k+2-k_\bc)p-1)\setminus\{jp-1:j\in \NN_1\}$ and let $X$ be a $\UMD$ Banach space. Let $\Delta_\bc$ be the Laplacian on $W^{k,p}_{\del,\bc}(\R^d_+, w_\gam;X)$ with domain
    \begin{equation*}
        D(\Delta_\bc):=W^{k+2,p}_{\del,\bc}(\R^d_+, w_\gam;X).
    \end{equation*}
    The following assertions hold.
    \begin{enumerate}[(i)]
        \item\label{it:thm1D1} For all $\lambda>0$, the operator $\lambda-\Delta_\bc$ is sectorial of angle $\om(\lambda-\del_\bc)=0$.
        \item\label{it:thm1D2} $(T^d_\bc(z))_{z\in \Sigma_\sigma}$ with $\sigma\in (0, \frac{\pi}{2})$ is an analytic $C_0$-semigroup generated by $\del_\bc$.
        \item\label{it:thm1calculus} For all $\lambda>0$, the operator $\lambda-\Delta_\bc$ has a bounded $\Hinf$-calculus of angle $\om_{\Hinf}(\lambda-\del_\bc)=0$.
    \end{enumerate}
    In addition, the semigroup $T^d_\bc(t)$ on $W^{k,p}_{\del, \bc}(\R^d_+,w_\gam;X)$ satisfies the following growth properties.
    \begin{enumerate}[resume*]
        \item\label{it:thm1D3} If $\gam< (2-k_\bc)p-1$, then $T^d_\bc(t)$ is bounded and assertions \ref{it:thm1D1} and \ref{it:thm1calculus} hold for $\lambda=0$ as well.
        \item\label{it:thm1D4} If $\gam>(2-k_\bc)p-1$, then $T^d_\bc(t)$ has polynomial growth and satisfies
        \begin{equation*}
            \|T^d_\bc(t)\|\eqsim 1+ t^{\alpha_{\bc,\gamma}-1},\quad  \alpha_{\bc,\gamma}:= \max\Big\{1, \frac{\gam+k_\bc p+1}{2p}\Big\},\,  t\geq 0.
        \end{equation*}
    \end{enumerate}
\end{theorem}
Theorem \ref{thm:intro_main_thm} is an extension of \cite[Theorems 1.1 \& 1.2]{LLRV24} which deal with the weight exponents
\begin{equation}\label{eq:gamtilde}
    \gam\in  \big((k-k_\bc)p-1, (k+2-k_\bc)p-1\big)\setminus\{(k+1-k_\bc)p-1\}.
\end{equation}
In particular, it was proved in \cite[Theorems 1.1(v) \& 1.2(v)]{LLRV24} that for any $\eps>0$ we have the upper estimate
        \begin{equation*}
            \|T^d_\bc(t)\|\lesssim 1+ t^{\alpha_{\bc, \gam}-1 + \eps},\quad t\geq 0.
        \end{equation*}
        Theorem \ref{thm:intro_main_thm} proves that this estimate actually holds with $\eps=0$ and the range of weight exponents \eqref{eq:gamtilde} is extended to $\gam>-1$. The optimal estimates for the semigroups in Theorem \ref{thm:intro_main_thm}\ref{it:thm1D4} are obtained in two steps. First, via explicit kernel representations we prove sharp $\lambda$-dependent estimates for the one-dimensional resolvent operator. Since resolvent estimates are equivalent to growth estimates on the semigroup (see Lemma \ref{lem:growth_semigroup_abstract}), this implies Theorem \ref{thm:intro_main_thm}\ref{it:thm1D4} for $d=1$. Secondly, using the factorisation of the semigroup in a normal and tangential part $T^d_\bc = T^{d-1}T_\bc^1$, we extend the one-dimensional result to arbitrary dimensions.
\begin{remark}\label{rem:intro}
    We make the following remarks about Theorem \ref{thm:intro_main_thm}.
    \begin{enumerate}[(i)]
        \item\label{it:rem:intro1} The range $\gam< (k+2-k_\bc)p-1$ is optimal in the sense that for larger $\gam$ the semigroups $ T_\bc^d$ are not well defined, see \cite[Example 6.9]{LLRV24}.
        \item\label{it:rem:intro2} The critical values $\gam=jp-1$ with $j\in \NN_1$ are excluded since in this case the domain characterisation of the semigroup generator fails. Nevertheless, for $j>2-k_\bc$ it is expected that the growth estimates on the semigroup remain true and one could identify the domain of the generator as a closure of a suitable space of test functions instead, cf. \cite[Remark 4.3]{LLRV24}.

        In the special case $\gam=(2-k_\bc)p-1$, we conjecture that the growth of the semigroup is logarithmic:
        \begin{equation*}
            \|T^d_\bc (t)\|\eqsim \big(1+ \log(1+t)\big)^\frac{1}{p'},\quad t\geq 0,
        \end{equation*}
        see Remark \ref{rem:critical-log} for a more elaborate discussion.
        \item\label{it:rem:intro3} The blow-up of the semigroup is actually caused by the growth of the weight $w_\gam(x)=|x_1|^\gam$ at infinity. In the case of bounded domains and the weight exponents \eqref{eq:gamtilde}, it was proved in \cite{LLRV25} that the semigroup is bounded. Similar techniques as developed in \cite{LLRV25} can be used to prove a version of Theorem \ref{thm:intro_main_thm} on bounded domains. 
        
        Alternatively, one could replace the weights $w_\gam$ in Theorem \ref{thm:intro_main_thm} by weights that behave like $|x_1|^\gam$ close to the boundary and that are constant at infinity. In this setting, we expect that the semigroup is also bounded independently of the smoothness index and weight exponent, see Remark \ref{rem:rho}.
        \item \emph{Necessity of the compatibility conditions.} If $-1 < \gam <(k-k_\bc)p-1$, then already on $\RR_+$ the realisation of the Laplacian on $W^{k,p}(\RR_+,w_\gam;X)$ without compatibility conditions is not sectorial, even after a positive shift; see Proposition \ref{prop:notsect}. For the unweighted setting $\gam=0$ this was already proved in \cite{DD11}, even for more general elliptic and mixed-order operators. Nevertheless, it should be noted that in general (for fixed $k$) \emph{less} compatibility conditions are needed if $\gamma$ is larger. Especially, for weights exponents \eqref{eq:gamtilde} (which is the setting as studied in \cite{LLRV24}), no compatibility conditions are needed at all.
    \end{enumerate}
\end{remark}

\subsubsection*{Outline} The outline of this paper is as follows. In Section \ref{sec:elliptic}
we first prove some density and elliptic regularity results. In addition, we give the proof of Theorem \ref{thm:intro_main_thm}\ref{it:thm1D1}, \ref{it:thm1D2} and \ref{it:thm1calculus}. In Section \ref{sec:1D} we prove optimal estimates for the one-dimensional resolvent operator on weighted Sobolev spaces. Finally, in Section \ref{sec:semigroup} we establish the optimal estimates on the heat semigroups and we complete the proof of Theorem \ref{thm:intro_main_thm}.

\subsection*{Notation}
We denote by $\NN_0$ and $\NN_1$ the sets of natural numbers starting at $0$ and $1$, respectively. 
For $a\in \R$, we use the notation $(a)_+:=\max\{a,0\}$. 
For $d\in\NN_1$ the half-space is given by $\RRdh=\RR_+\times\RR^{d-1}$, where $\RR_+=(0,\infty)$ and for $x\in \RRdh$ we write $x=(x_1,\tilde{x})$ with $x_1\in \RR_+$ and $\tilde{x}\in \RR^{d-1}$. 

For two topological vector spaces $X$ and $Y$, the space of continuous linear operators is denoted by $\mc{L}(X,Y)$ and $\mc{L}(X):=\mc{L}(X,X)$. Unless specified otherwise, $X$ will always denote a Banach space with norm $\|\cdot\|_X$ and the dual space is $X':=\mc{L}(X,\CC)$.

For a linear operator $A:X\supseteq D(A)\to X$ on a Banach space $X$, we denote by $\sigma(A)$ and $\rho(A)$  the spectrum and resolvent set, respectively. For $\lambda\in\rho(A)$, the resolvent operator is given by $R(\lambda,A)=(\lambda-A)^{-1}\in \mc{L}(X)$.

We write $f\lesssim g$ (resp. $f\gtrsim g$) if there exists a constant $C>0$, possibly depending on parameters which will be clear from the context or will be specified in the text, such that $f\leq Cg$ (resp. $f\geq Cg$). Furthermore, $f\eqsim g$ means $f\lesssim g$ and $g\lesssim  f$.\\

For an open and non-empty $\OO\subseteq \RR^d$ and $\ell\in\NN_0\cup\{\infty\}$, the space $C^\ell(\OO;X)$ denotes the space of $\ell$-times continuously differentiable functions from $\OO$ to some Banach space $X$. As usual, this space is equipped with the compact-open topology. 

Let $\Cc^{\infty}(\OO;X)$ be the space of compactly supported smooth functions on $\OO$ equipped with its usual inductive limit topology. The space of $X$-valued distributions is given by $\mc{D}'(\OO;X):=\mc{L}(\Cc^{\infty}(\OO);X)$. Moreover, $\Cc^{\infty}(\overline{\OO};X)$ is the space of smooth functions with their support in a compact set contained in $\overline{\OO}$.

We denote the Schwartz space by $\SS(\RRd;X)$, and $\SS'(\RRd;X):=\mc{L}(\SS(\RRd);X)$ is the space of $X$-valued tempered distributions. For $f\in \SS(\RRd;X)$ we define the $d$-dimensional Fourier transform $(\mc{F} f)(\xi):=\hat{f}(\xi):=\int_{\RRd} f(x)e^{-\ii x\cdot\xi}\dd x$ for $\xi\in \RR^d$, which extends to $\SS'(\RRd;X)$ by duality.

\section{Sectoriality and functional calculus}\label{sec:elliptic}

In this section, we prove Theorem \ref{thm:intro_main_thm}\ref{it:thm1D1}, \ref{it:thm1D2} and \ref{it:thm1calculus} which concern sectoriality and boundedness of the functional calculus. We first prove some density results in Section \ref{subsec:density} and focus on elliptic regularity in Section \ref{subsec:elliptic}. Finally, the proof of Theorem \ref{thm:intro_main_thm}\ref{it:thm1D1}, \ref{it:thm1D2} and \ref{it:thm1calculus} is given in Section \ref{subsec:calculus}.\\

For completeness, we first give the precise definitions of sectoriality and boundedness of the $\Hinf$-calculus. For a more elaborate discussion of these notions, we refer to \cite{Ha06} and \cite[Chapter 10]{HNVW17}.
\begin{definition}
    An injective, closed linear operator $(A, D(A))$ with dense domain and dense range on a Banach space $X$ is called \emph{sectorial} if there exists an $\om\in (0, \pi)$ such that $\sigma(A)\subseteq \overline{\Sigma_{\om}}$ and 
\begin{equation*}
    \sup_{\lambda\in \CC\setminus\overline{\Sigma_\om}}\|\lambda R(\lambda, A)\|<\infty.
\end{equation*}
The angle of sectoriality $\om(A)$ is the infimum over all possible $\om$.
\end{definition}
Resolvent estimates are related to growth of the corresponding semigroup. 
In particular, if an operator is sectorial (of angle $<\frac \pi 2$), then the negative operator generates a bounded analytic semigroup. We recall the following extension of this result, which relates the blow-up of the resolvent estimate to the precise growth of the semigroup.

\begin{lemma}[{\cite[Lemma 2.3]{LLRV24}}]\label{lem:growth_semigroup_abstract}
  Let $A$ be a linear operator on a Banach space $X$ and let $\alpha\geq 1$. The following are equivalent.
  \begin{enumerate}[(i)]
    \item \label{it:lem:growth1} There exist $\om\in(0,\frac{\pi}{2})$ and $C_1>0$ such that $\sigma(A)\subseteq \overline{\Sigma_{\om}}$ and
         \begin{equation*}
            |\lambda|\|(\lambda + A)^{-1}\|_{\mc{L}(X)}\leq C_1 (1+ |\lambda|^{1-\alpha}),\qquad \lambda\in \Sigma_{\pi-\om}.
         \end{equation*}
    \item\label{it:lem:growth2} There exist $\eta\in(0,\frac{\pi}{2})$ and $C_2>0$ such that $-A$ generates an analytic $C_0$-semigroup on $\Sigma_{\eta}$ and
    \begin{equation*}
        \|e^{-zA}\|_{\mc{L}(X)}\leq C_2(1+|z|^{\alpha-1}),\qquad z\in \Sigma_{\eta}.
    \end{equation*}
  \end{enumerate}
\end{lemma}

We now introduce the functional calculus. Let $\om\in(0,\pi)$, then $H^1(\Sigma_{\om})$ is the Hardy space of all holomorphic functions $f:\Sigma_{\om}\to \CC$ such that
\begin{equation*}
  \|f\|_{H^1(\Sigma_{\om})}:=\sup_{|\nu|<\om}\|t\mapsto f(e^{\ii\nu}t)\|_{L^1(\RR_+,\frac{\mathrm{d}t}{t})}<\infty.
\end{equation*}
Moreover, let $H^{\infty}(\Sigma_{\om})$ be the Hardy space of all bounded holomorphic functions on the sector with norm
\begin{equation*}
  \|f\|_{H^{\infty}(\Sigma_{\om})}:=\sup_{z\in\Sigma_{\om}}|f(z)|.
\end{equation*}

\begin{definition}
  Let $A$ be a sectorial operator on a Banach space $X$ and let $\om\in(\om(A),\pi)$, $\nu\in(\om(A), \om)$ and $f\in H^{1}(\Sigma_{\om})$. We define the operator
  \begin{equation*}
    f(A):=\frac{1}{2\pi \ii}\int_{\partial \Sigma_{\nu}}f(z)R(z,A)\dd z,
  \end{equation*}
  where $\d \Sigma_{\nu}$ is traversed downwards. The operator $A$ has a \emph{bounded $\Hinf(\Sigma_{\om})$-calculus} if there exists a $C>0$ such that
  \begin{equation*}
    \|f(A)\|\leq C\|f\|_{\Hinf(\Sigma_{\om})}\quad \text{ for all }f\in H^1(\Sigma_{\om})\cap\Hinf(\Sigma_{\om}).
  \end{equation*}
  The angle of the $\Hinf$-calculus $\om_{\Hinf}(A)$ is defined as the infimum over all possible $\om>\om(A)$.
\end{definition}

\subsection{Density results}\label{subsec:density}
We start by recalling the following from \cite[Section 3]{LLRV24}. 
For $j\in\NN_0$ and $X$ a Banach space, define
\begin{equation*}
  C^{\infty}_{{\rm c},j}(\overline{\RRdh};X):=\{f\in \Cc^{\infty}(\overline{\RRdh};X): \d_1^j f\in \Cc^{\infty}(\RRdh;X)\}.
\end{equation*}
The condition $\d_1^j f\in \Cc^{\infty}(\RRdh;X)$ implies that $\d^{\alpha} f\in \Cc^{\infty}(\RRdh;X)$ for all $\alpha=(\alpha_1,\tilde{\alpha})\in\NN_0\times \NN_0^{d-1}$ with $\alpha_1\geq j$.

\begin{lemma}[{\cite[Lemma 3.4]{LLRV24}}]\label{lem:densityLLRV}
  Let $p\in (1,\infty)$, $j,k\in \NN_0$ such that $k\geq j$ and  $\gam>(k-j)p-1$, and let $X$ be a Banach space. Then $C^{\infty}_{{\rm c},j}(\overline{\RRdh};X)$ is dense in $W^{k,p}(\RRdh,w_{\gam};X)$.
\end{lemma}
In particular, Lemma \ref{lem:densityLLRV} implies that the  space of test functions $\Cc^\infty(\RRdh;X)= C^{\infty}_{{\rm c},0}(\overline{\RRdh};X)$ is dense in $W^{k,p}(\RRdh, w_\gam;X)$ if $\gam>kp-1$.\\

We introduce a similar class of test functions to deal with the compatibility conditions required for the sectoriality of the Laplacian on weighted Sobolev spaces.
Let $j\in \NN_0$, $\bc\in \{\Dir, \Neu\}$ and let $X$ be a Banach space. We define the space of test functions given by
\begin{equation*}
    C^{\infty}_{\del, \bc, j}(\overline{\RRdh};X):=\{f\in C^\infty_{{\rm c}, j}(\overline{\RRdh};X): (\d_1^{k_\bc+2r}f)|_{\d\RRdh}=0, \, r\in \NN_0,\, k_\bc + 2r < j  \},
\end{equation*}
where $k_\Dir:=0$ and $k_\Neu:=1$.
For notational convenience, we also set $C^{\infty}_{\del, \bc, -1}(\overline{\RRdh};X)=C^\infty_{{\rm c}, -1}(\overline{\RRdh};X)=\Cc^\infty(\RRdh;X)$. 

We consider the setting in which this space of test functions is dense in $W^{k,p}_{\del,\bc}(\RRdh, w_\gam;X)$. To this end, for $p\in (1,\infty)$, $k\in \NN_0$ and $\gam\in (-1, (k+2-k_\bc)p-1)\setminus\{jp-1:j\in \NN_1\}$, we define $j_*\in \NN_0\cup \{-1\}$ by
\begin{equation}\label{eq:jstar}
    j_*:= \Big\lceil k- \frac{\gam+1}{p}\Big\rceil.
\end{equation}
Then $j_*$ is the unique integer such that $(k-j_*)p-1<\gam<(k-j_*+1)p-1$. In particular, in $W^{k,p}(\RRdh, w_\gam;X)$ all traces of order $r<j_*$ exist, while no trace operators of order $r\geq j_*$ are continuous. We have the following density result.

\begin{lemma}\label{lem:density}
    Let $p\in (1,\infty)$, $k\in \NN_0$, $\bc\in \{\Dir, \Neu\}$, $\gam\in (-1, (k+2-k_\bc)p-1)\setminus\{jp-1:j\in \NN_1\}$ and let $X$ be a Banach space. Let $j_*$ be as defined in \eqref{eq:jstar}. Then $C^{\infty}_{\del, \bc, j_*}(\overline{\RRdh};X)$ is dense in $W^{k,p}_{\del, \bc}(\RRdh, w_\gam;X)$.
\end{lemma}
\begin{proof} If $j_*\leq 0$, then $\gam>kp-1$ and the result follows directly from Lemma \ref{lem:densityLLRV}. We may thus assume that $j_*\geq 1$.
   Let $f\in W^{k,p}_{\del, \bc}(\RRdh, w_\gam;X)$ and $\eps>0$. By \cite[Theorem 7.2 \& Remark 11.12(iii)]{Ku85}, which also holds in the vector-valued case, and a standard cut-off argument, there exists a $g\in \Cc^\infty(\overline{\RRdh};X)$ such that
   \begin{equation}\label{eq:dens1}
       \|f-g\|_{W^{k,p}(\RRdh, w_\gam;X)}<\eps.
   \end{equation}
   Define the index set $I_\bc:=\{k_\bc+2r: r\in\NN_0,\, k_\bc+2r<j_* \}$ and the (system of) trace operators $\mc{B}_{\bc}:=(\Tr_m)_{m\in I_{\bc}}$. By \cite[Theorem~5.6]{Ro25}, we have that
\begin{equation*}
    \mc{B}_{\bc}: W^{k,p}(\RRdh, w_\gam;X)\to \prod_{m\in I_\bc} B_{p,p}^{k-m-\frac{\gam+1}{p}}(\RR^{d-1};X)
\end{equation*}
is continuous and surjective. Let $\ext_{\bc}$ denote the corresponding right inverse to $\mc{B}_{\bc}$. Repeating the proof of \cite[Proposition 4.12]{Ro25} using the approximation $h:= g - \ext_{\bc}(\mc{B}_\bc g)$ (or $h:=g$ if $I_\bc=\varnothing$) and \eqref{eq:dens1}, we obtain that $\|f-h\|_{W^{k,p}(\RRdh, w_\gam;X)}< C\eps$. Note that $h\in C^\infty(\overline{\RRdh};X)$ with $(\d_1^m h)|_{\d\RRdh}=0$ for $m\in I_\bc$, and after another cut-off argument, we find that 
\begin{equation*}
    \mc{T}:=\big\{h\in \Cc^\infty (\overline{\RRdh};X): (\d_1^{k_\bc+2r}h)|_{\d\RRdh}=0, k_\bc + 2r < j_* \big\} \text{ is dense in }W^{k,p}_{\del, \bc}(\RRdh, w_\gam;X).
\end{equation*}
It remains to approximate functions in $\mc{T}$ with functions that are additionally in $C^\infty_{{\rm c}, j_*}(\overline{\RRdh};X)$. Let $f\in \mc{T}$. In addition, let $\phi\in C^\infty(\RR_+)$ be such that $\phi=0$ on $[0, \half]$ and $\phi =1 $ on $[1,\infty)$ and set $\phi_n(x_1):=\phi(nx_1)$. Define
\begin{equation*}
     P(x_1,\widetilde x)
      :=\sum_{\ell=0}^{j_*-1}
        \frac{x_1^\ell}{\ell!}\Tr_\ell f(\widetilde x)\quad \text{ and }\quad f_n:=P+\phi_n(f-P).
\end{equation*}
Then $f_n=f$ for $x_1\geq n^{-1}$, while $f_n=P$ near the
boundary. Consequently,
\[
    \partial_1^{j_*}f_n=0
    \quad\text{near }\partial\RRdh
    \quad\text{ and }\quad
    \Tr_\ell f_n=\Tr_\ell f
    \quad(0\leq\ell<j_*).
\]
Thus $f_n\in C^\infty_{\del,\bc,j_*}(\overline{\RRdh};X)$ and by the product rule we obtain
\[
    \|f_n-f\|_{W^{k,p}(\RRdh,w_\gam;X)}
    \lesssim
    n^{-(j_*-k+\frac{\gam+1}{p})},
\]
which can be proved with similar estimates as in the proof of \cite[Proposition 4.13]{Ro25}.
Since  \(j_*>k-(\gam+1)/p\), we have  $\|f_n-f\|_{W^{k,p}(\RRdh,w_\gam;X)}\to 0$ as $n\to \infty$, which completes the proof.
\end{proof}

\subsection{Elliptic regularity}\label{subsec:elliptic}
We first prove solvability of the resolvent equation associated to the Laplacian with smooth right-hand side satisfying compatibility conditions. This can be proved using a suitable extension from $\RRdh$ to $\RR^d$. 
We define the odd (for Dirichlet boundary conditions $\bc=\Dir$) and even (for Neumann boundary conditions $\bc=\Neu$) extension in $x_1=0$ by
\begin{equation*}
    (\mc{E}_\bc f)(x_1,\tilde{x}) = \begin{cases}
        f(x_1, \tilde{x})&\mbox{ if }x_1\geq 0,\\
        (-1)^{k_\bc+1}f(-x_1, \tilde{x})&\mbox{ if }x_1<0.
    \end{cases}
\end{equation*}
Recall that $k_\Dir:=0$ and $k_\Neu:=1$. Moreover, the Schwartz space on $\RRdh$ is defined as $\mc{S}(\RRdh;X):=\{u|_{\RRdh}: u\in \mc{S}(\RRd;X)\}$. 
\begin{lemma}\label{lem:smoothRHS}
    Let $\bc\in \{\Dir, \Neu\}$, $j\in \NN_1$, $\om\in (0,\pi)$ and let $X$ be a $\UMD$ Banach space. Then for every $f\in C^\infty_{\del, \bc, j}(\overline{\R^d_+};X)$ and $\lambda \in \Sigma_{\pi-\om}$ there exists a unique $u \in \mc{S}(\RRdh;X)$ such that
    \begin{equation*}
        \lambda u -\del u =f,\qquad (\del^{r} \d_1^{k_\bc} u)(0, \cdot)=0 \text{ for all }r\in \NN_0.
    \end{equation*}
    In particular, the solution is given by
    \begin{equation*}
    u =  \Big(\mc{F}^{-1}\Big[\xi \mapsto \frac{ (\mc{F}\mc{E}_\bc f)(\xi)}{\lambda+ |\xi|^2}\Big]\Big)\Big|_{\RRdh}.
\end{equation*}
\end{lemma}
\begin{proof}
    Since $\d_1^j f\in \Cc^\infty(\RRdh;X)$ and $(\d_1^{k_\bc+2r}f)(0,\cdot)=0$, there exists a $\delta>0$ such that
    \begin{equation*}
        f(x_1, \tilde{x})=\sum_{\substack{0\leq m< j\\ m+k_\bc\text{ odd}}} \frac{x_1^m}{m!}(\d_1^m f)(0, \tilde{x}),\qquad x_1\in [0,\delta),\, \tilde{x}\in \RR^{d-1}.
    \end{equation*}
    Applying the odd or even extension to $\RR^d$, we find that $\mc{E}_\bc f$ is given by the same formula, i.e.,
        \begin{equation*}
        (\mc{E}_\bc f)(x_1, \tilde{x})=\sum_{\substack{0\leq m< j\\ m+k_\bc\text{ odd}}} \frac{x_1^m}{m!}(\d_1^m f)(0, \tilde{x}),\qquad |x_1|<\delta,\, \tilde{x}\in \RR^{d-1}.
    \end{equation*}
    Therefore, $\mc{E}_\bc f$ is also smooth in $x_1=0$ and we find $\mc{E}_\bc f \in \Cc^\infty(\RRd;X)\subseteq \mc{S}(\RR^d;X)$. Taking the Fourier transform of the resolvent equation $\lambda u -\del u =\mc{E}_\bc f $ on $\RR^d$ yields a solution $u_\bc \in \mc{S}(\RR^d;X)$ given by
    \begin{equation*}
    u_\bc =  \mc{F}^{-1}\Big[\xi \mapsto \frac{ (\mc{F}\mc{E}_\bc f)(\xi)}{\lambda+ |\xi|^2}\Big].
\end{equation*}
Now $\d_1^{k_\bc} u_\bc$ is odd and therefore $\del^r \d_1^{k_\bc} u_\bc$ is odd as well. This implies that $(\del^r \d_1^{k_\bc}u_\bc )(0, \cdot)=0$. By setting $u:= u_\bc |_{\RRdh}$, we obtain a solution to the resolvent equation on $\RRdh$. The uniqueness follows from \cite[Corollary 4.7]{LLRV24} and \cite[Corollary 4.8]{LLRV24} for Dirichlet and Neumann boundary conditions, respectively.
\end{proof}

To continue, we combine Lemma \ref{lem:smoothRHS} and density from Lemma \ref{lem:density} to prove elliptic regularity for the Laplacian on weighted Sobolev spaces with compatibility conditions, which extends \cite[Propositions 5.4 \& 5.6]{LLRV24}. It should be noted that  in the elliptic regularity estimate below, the constant $C$ depends on $\lambda$ and blows up as $\lambda\downarrow 0$. In the subsequent section, we will establish the precise growth of the resolvent in terms of $\lambda$. 
\begin{proposition}\label{prop:sect_nonoptimal}
    Let $p\in (1,\infty)$, $k\in \NN_0$, $\bc\in \{\Dir, \Neu\}$, $\gam\in (-1, (k+2-k_\bc)p-1)\setminus\{jp-1:j\in \NN_1\}$, $\omega\in (0,\pi)$, $\lambda_0>0$ and let $X$ be a $\UMD$ Banach space. Then for all $f\in W^{k,p}_{\del, \bc}(\RRdh, w_{\gam};X)$ and $\lambda\in \Sigma_{\pi-\om}$ with $|\lambda|>\lambda_0$ there exists a unique $u\in W^{k+2,p}_{\del, \bc}(\RRdh, w_{\gam};X)$ to $\lambda u -\del_\bc u =f$. Moreover, this solution satisfies
    \begin{equation*}
         \sum_{|\beta|\leq2}
        |\lambda|^{1-\frac{|\beta|}{2}}\|\d^\beta u\|_{W^{k,p}(\RRdh, w_{\gam};X)}\leq C\|f\|_{W^{k,p}(\RRdh, w_\gam;X)},
    \end{equation*}
    where the constant $C>0$ only depends on $p,k,\gam, \omega, \bc, \lambda_0, d$ and $X$.
\end{proposition}
\begin{proof} For $\gam>kp-1$ (i.e., $ j_*\leq 0$), the result is already contained in \cite[Proposition 5.4 \& 5.6]{LLRV24}. Let $j_*\geq 1$, we will differentiate the resolvent equation to reduce to the case $j_*=0$ (see also \cite[Lemmas 4.5 \& 4.10]{LLRV25}).
    Take $f\in C_{ \del,\bc, j_*}^\infty(\overline{\RRdh};X)$ and let $u\in \mc{S}(\RRdh;X)$ be the solution to the resolvent equation from Lemma \ref{lem:smoothRHS}, which satisfies $\Tr \d_1^{k_\bc+2r}u=0$ for all $r\in \NN_0$. Set $\ell:= k- j_*$, i.e., $\gam$ satisfies $\ell p-1< \gam< (\ell+1)p-1$. Let $\alpha=(\alpha_1,\tilde{\alpha})\in \NN_0\times \NN_0^{d-1}$ be such that $|\alpha|\leq k-\ell$. Then $v_\alpha:=\d^\alpha u $ satisfies the equation $\lambda v_\alpha-\del v_\alpha = \d^\alpha f$ with one of the following boundary conditions
    \begin{equation*}
        \begin{cases}
            v_\alpha(0,\cdot)=0& \mbox{ if }k_\bc + \alpha_1\text{ is even,}\\
            (\d_1 v_\alpha)(0,\cdot)=0& \mbox{ if }k_\bc+\alpha_1\text{ is odd.}\\
        \end{cases}
    \end{equation*}
    We may therefore apply \cite[Proposition~5.4]{LLRV24} with parameters
$(\ell,\gam-\ell p)$ if $v_\alpha$ satisfies a Dirichlet boundary condition, and
\cite[Proposition~5.6]{LLRV24} with parameters
$(\ell-1,\gam-(\ell-1)p)$ if $v_\alpha$ satisfies a Neumann boundary condition. In both cases, this yields the estimates
    \begin{equation*}
        \sum_{|\beta|\leq 2}|\lambda|^{1-\frac{|\beta|}{2}}\|\d^\beta v_\alpha\|_{W^{\ell,p}(\RRdh,w_\gam;X)}
        \leq C
        \|\d^\alpha f\|_{W^{\ell,p}(\RRdh,w_\gam;X)},
        \qquad |\alpha|\leq k-\ell.
    \end{equation*}
    Consequently, 
    \begin{align*}
        \sum_{|\beta|\leq 2}|\lambda|^{1-\frac{|\beta|}{2}}\|\d^\beta u\|_{W^{k,p}(\RRdh,w_\gam;X)}
        &\lesssim
        \sum_{|\beta|\leq 2}\sum_{|\alpha|\leq k-\ell}|\lambda|^{1-\frac{|\beta|}{2}}
        \|\d^\beta v_\alpha\|_{W^{\ell,p}(\RRdh,w_\gam;X)}
        \\
        &\lesssim
        \sum_{|\alpha|\leq k-\ell}
        \|\d^\alpha f\|_{W^{\ell,p}(\RRdh,w_\gam;X)}
        \lesssim
        \|f\|_{W^{k,p}(\RRdh,w_\gam;X)}.
    \end{align*}
    A density argument using Lemma \ref{lem:density}, similar to the proof of \cite[Proposition 5.4]{LLRV24}, yields the desired result. In addition, note that the uniqueness also follows from \cite[Proposition 5.4 \& 5.6]{LLRV24}.
\end{proof}

Using the elliptic regularity from Proposition \ref{prop:sect_nonoptimal}, we can characterise the domain of powers of $\lambda-\del_\bc$. In the proposition below, we also need a weighted Sobolev space with negative smoothness. For $\gam\in (-1,p-1)$ and $X$ a reflexive Banach space, we define
\begin{equation*}
    W^{-1, p}_{\del,\Dir}(\RRdh, w_\gam;X):= W^{-1,p}(\RRdh, w_\gam;X):= \big(W^{1,p'}(\RRdh, w_{\gam'};X')\big)',
\end{equation*}
where $\gam':=\frac{-\gam}{p-1}\in (-1, p'-1)$ is the dual weight. This definition can be shown to be equivalent to the definition in \cite[Section 4.3]{LLRV24}.
\begin{proposition}\label{prop:domains-integer-powers}
    Let $p\in(1,\infty)$, $k\in\NN_0$, $\bc\in\{\Dir,\Neu\}$, $\gam\in
   (-1,(k+2-k_\bc)p-1)
   \setminus\{jp-1:j\in\NN_1\}$ and let $X$ be a $\UMD$ Banach space. Let $\del_\bc$ be the Laplacian on $W^{k,p}_{\del,\bc}(\RRdh, w_\gam;X)$ with domain
   \begin{equation*}
       D(\del_\bc):= W^{k+2,p}_{\del,\bc}(\RRdh, w_\gam;X).
   \end{equation*}
   Let $\lambda>0$ and consider $A:=\lambda-\del_\bc$ on $W^{k,p}_{\del,\bc}(\RRdh, w_\gam;X)$. Then, for every $\ell\in\NN_0$,
\begin{equation*}
   D(A^\ell)=W_{\del,\bc}^{k+2\ell,p}(\RRdh,w_\gam;X)\quad \text{ with equivalent norms}.
\end{equation*} 
Additionally, if $k=-1$, $\bc=\Dir$ and $\gam\in (-1, p-1)$, then $A:=\lambda-\del_\Dir$ on $W^{-1,p}_{\del,\Dir}(\RRdh, w_\gam;X)$ satisfies
\begin{equation*}
    D(A^\ell)=W_{\del,\Dir}^{-1+2\ell,p}(\RRdh,w_\gam;X)\quad \text{ with equivalent norms}.
\end{equation*}
\end{proposition}
\begin{proof} We only consider the case $k\geq 0$, while the case $k=-1$ can be proved similarly using the elliptic regularity results from \cite[Proposition 4.9]{LLRV24}.
    For notational convenience, we define $Y_\ell:=W^{k+2\ell, p}_{\del,\bc}(\RRdh, w_\gam;X)$. Consider the operators
\[
  A_\ell:Y_{\ell+1}\to Y_\ell,
  \qquad A_\ell u:=\lambda u-\del_\bc u,
\]
which are consistent. To see that $A_\ell$ maps into $Y_\ell$, note that by properties of the trace we have
\begin{equation*}
    \Tr\d_1^{k_\bc+2r}(\lambda u -\del_\bc u) =0,\qquad u\in Y_{\ell+1},
\end{equation*}
for any $r\in\NN_0$ such that $k_\bc + 2r< k+2\ell-\frac{\gam+1}{p}$.
By Proposition \ref{prop:sect_nonoptimal} (applied to $k+2\ell$ instead
of $k$) we have that
\begin{equation}\label{eq:Aj-isomorphism}
  A_\ell:Y_{\ell+1}\to Y_\ell
  \quad\text{is an isomorphism}.
\end{equation}
To prove that $D(A^\ell)=Y_\ell$ we argue by induction on $\ell\geq 0$. Note that $D(A^0)=Y_0$ and $D(A)=Y_1$ by definition. Suppose that $D(A^\ell)=Y_\ell$ for some $\ell\geq 0$. The induction hypothesis implies that 
\[
  D(A^{\ell+1})
  =
  \{u\in D(A): Au\in D(A^\ell)\} = \{u\in Y_1: A_0 u\in Y_\ell\}.
\]
The embedding $Y_{\ell+1}\hookrightarrow D(A^{\ell+1})$ is now straightforward to check using
$A_\ell Y_{\ell+1}\hookrightarrow Y_\ell$ and consistency.

Conversely, let $u\in Y_1$ and $f:=A_0u\in Y_\ell$. By
\eqref{eq:Aj-isomorphism}, there is a unique
$v\in Y_{\ell+1}$ such that $A_\ell v=f$. Consistency gives
$A_0v=f=A_0u$, and injectivity of $A_0$ yields $u=v$. Hence
$u\in Y_{\ell+1}$. This completes the induction.
\end{proof}

\subsection{The proof of Theorem \ref{thm:intro_main_thm}\ref{it:thm1D1}-\ref{it:thm1calculus}}\label{subsec:calculus}
To conclude this section, we give the proof of Theorem \ref{thm:intro_main_thm}\ref{it:thm1D1}, \ref{it:thm1D2} and \ref{it:thm1calculus}.

\begin{proof}[Proof of Theorem \ref{thm:intro_main_thm}\ref{it:thm1D1}, \ref{it:thm1D2} and \ref{it:thm1calculus}] We prove the three statements separately.

\textit{Proof of \ref{it:thm1D1}. }
Arguing as in \cite[Section 5.3]{LLRV24}, using
Proposition \ref{prop:sect_nonoptimal}, we find that
$\lambda-\del_\bc$ is sectorial of angle zero for every $\lambda>0$.

\textit{Proof of \ref{it:thm1D2}. }As a consequence of \ref{it:thm1D1} and a shift, $\del_\bc$ generates an analytic $C_0$-semigroup on
$W^{k,p}_{\del,\bc}(\RRdh,w_\gam;X)$ (see Lemma \ref{lem:growth_semigroup_abstract}). It remains to identify this
semigroup with $T^d_\bc$. If $j_*\leq k_\bc$, this follows directly from
\cite[Section 6.1]{LLRV24}, with the parameter choices used in the
proof of Proposition \ref{prop:sect_nonoptimal}. Let $j_*$ be as in \eqref{eq:jstar}. If $j_*> k_\bc$, set
$\ell:=k-j_*$ and consider the lower-order realisation of $\del_\bc$
on $W^{\ell,p}(\RRdh,w_\gam;X)$. By uniqueness in Proposition
\ref{prop:sect_nonoptimal}, its resolvent is consistent with the
resolvent on $W^{k,p}_{\del,\bc}(\RRdh,w_\gam;X)$ as can be proved similarly to \cite[Lemmas 6.4 \& 6.5]{LLRV24}. Arguing as in the proof of the generator identification in \cite[Theorems 6.1 \& 6.2]{LLRV24}, shows that the semigroup generated by $\del_\bc$ coincides with $T^d_\bc$. 

\textit{Proof of \ref{it:thm1calculus}. }Let $j_*$ be as in \eqref{eq:jstar} and set $\ell:=k-j_*$, i.e., $\ell p -1 < \gam < (\ell+1)p-1$. If $j_*\leq k_\bc$, then $W_{\del,\bc}^{k,p}(\RRdh,w_\gam;X)
    =W^{k,p}(\RRdh,w_\gam;X)$ and $W_{\del,\bc}^{k+2,p}(\RRdh,w_\gam;X)
    =W_{\bc}^{k+2,p}(\RRdh,w_\gam;X)$. It follows from \cite[Theorems 1.1 \& 1.2]{LLRV24} (for $\bc=\Dir$ and $\bc=\Neu$, respectively) that $A:=\lambda-\del_\bc$ has a bounded $\Hinf$-calculus of angle zero. 

    Assume now that $j_*>k_\bc$ and set
\[
    N:=\left\lfloor\frac{j_*+1-k_\bc}{2}\right\rfloor\in\NN_1.
\]
Then $k-k_\bc-2N=\ell+j_*-k_\bc-2N\in \{\ell-1, \ell\}$ and thus
\begin{equation*}
    \gam - (k-k_\bc-2N)p \in (-1, 2p-1)\setminus\{p-1\}.
\end{equation*}
Consider the realisation $B:=\lambda-\del_\bc$ on $W^{k-2N,p}_{\del,\bc}(\RRdh, w_\gam;X)$ with domain $D(B)=W^{k-2N+2,p}_{\del,\bc}(\RRdh, w_\gam;X)$. Note that $k-2N\in \NN_0\cup\{-1\}$ if $\bc=\Dir$, and $k-2N\in \NN_0$ if $\bc=\Neu$. Due to the choice of $N$, there are no compatibility conditions in the space and again by \cite[Theorems 1.1 \& 1.2]{LLRV24} we find that $B$ has a bounded $\Hinf$-calculus of angle zero. By Proposition~\ref{prop:domains-integer-powers}, we find
\begin{equation}\label{eq:Hinf-power-domains}
    D(B^N)
    =W_{\del,\bc}^{k,p}(\RRdh,w_\gam;X)\quad \text{ and }\quad 
    D(B^{N+1})
    =W_{\del,\bc}^{k+2,p}(\RRdh,w_\gam;X).
\end{equation}
Since $\lambda>0$, $B$ is
invertible and 
\[
    J:=B^N:
    D(B^N) \to 
    W_{\del,\bc}^{k-2N,p}(\RRdh,w_\gam;X)
\]
is an isomorphism. By \eqref{eq:Hinf-power-domains}, $A$ is the part of $B$ in
$D(B^N)$ and 
\[
    D(A)=D(B^{N+1}) =\{u\in D(B^N):Bu\in D(B^N)\},
\]
and $Au=Bu$ on this domain. Moreover,
\(
    J(D(A))=D(B)
\)
and, for $u\in D(A)$,
\(
    JAu=B^NBu=BB^Nu=BJu.
\)
Hence, we obtain
\(
    A=J^{-1}BJ,
\)
which implies that 
\[
    R(z,A)=J^{-1}R(z,B)J,
    \qquad z\in\rho(B).
\]
Let $\om \in (0, \pi)$, the definition of the functional calculus gives
\[
    f(A)=J^{-1}f(B)J,
    \qquad f\in H^1(\Sigma_\om)\cap H^\infty(\Sigma_\om).
\]
Therefore,
\[
 \|f(A)\|_{
   \mc L(W_{\del,\bc}^{k,p}(\RRdh,w_\gam;X))}
 \leq
 \|J^{-1}\|\,
 \|f(B)\|_{
   \mc L(W_{\del,\bc}^{k-2N,p}(\RRdh,w_\gam;X))}
 \,\|J\|
 \lesssim_\om
 \|f\|_{H^\infty(\Sigma_\om)}.
\]
Since $\om\in(0,\pi)$ was arbitrary, this proves that the operator $A$ has a
bounded $\Hinf$-calculus of angle zero.
\end{proof}

\section{Optimal resolvent estimates in one dimension}\label{sec:1D}
In this section, we derive sharp estimates on the resolvent of the one-dimensional Dirichlet and Neumann Laplace operators. \\

We first introduce the kernel representation of the one-dimensional resolvent operators. 
Let $\bc\in \{\Dir,\Neu\}$ and recall that $k_\Dir:=0$ and $k_\Neu:=1$. For $m\in \NN_0$, we define the kernel
\begin{equation*}
    E_m(\lambda,x,y):=\frac{1}{2\sqrt{\lambda}}\Big(e^{-\sqrt{\lambda}|x-y|}+ (-1)^{m+1}e^{-\sqrt{\lambda}|x+y|}\Big),\quad \lambda\in \Sigma_{\pi},\, x, y>0, 
\end{equation*}
where we use the principal branch of the complex logarithm.
Note that $E_{m}=E_{m\,(\mod 2)}$ and we will write $E_\bc:= E_{k_\bc}$. In addition, for later reference, this kernel satisfies the following elementary properties:
\begin{equation}\label{eq:prop1DkernelR}
    \begin{aligned}
        \d_x E_m (\lambda, x,y)&= -\partial_{y}E_{m+1}(\lambda, x,y)\quad \text{ almost everywhere},\\
        E_{2m}(\lambda, x,0)&=0\quad \text{ and }\quad E_{2m+1}(\lambda,x,0)=\frac{1}{\sqrt{\lambda}}e^{-\sqrt{\lambda}x}.
    \end{aligned}
\end{equation}
For $m\in \NN_0$ and $f\in \Cc^\infty([0,\infty);X)$ we define the integral operator
\begin{equation*}
    R_{m}(\lambda)f(x)=\int_0^\infty E_{m}(\lambda, x,y)f(y) \dd y,\qquad x>0.
\end{equation*}
Again, we write $R_\bc:= R_{k_\bc}$. The integral operator $R_\bc$ coincides with the resolvent operator of the Laplacian $\Delta_\bc$ and extends to a bounded operator on certain weighted Sobolev spaces as we will prove in this section. To be precise, in Section \ref{subsec:1Dresovent_upper} we first prove sharp estimates for the resolvents $R_\bc$ on spaces with weights $w_{\gam}$ with $\gam=\tgam+kp$ and $\tgam\in (-1, 2p-1)\setminus\{p-1\}$. This setting is already considered in \cite{LLRV24}. In Section \ref{subsec:1Dresovent_lower} we extend the resolvent estimates to weights $w_\gam$ with more general weight exponents.

\subsection{The case of large weight exponents}\label{subsec:1Dresovent_upper}
For $m\in \NN_0$ and $\tilde{\gam}\in \R$ define
\begin{equation}\label{eq:def_hlambda}
    h_{m,\tgam}(\lambda):=1+ |\lambda|^{-\frac{(\tgam+mp-2p+1)_+}{2p}}.
\end{equation}
Recall that $(a)_+:=\max\{a,0\}$ for any $a\in \RR$.
The following theorems prove sharp resolvent estimates for the Dirichlet and Neumann Laplacian in one dimension. 
\begin{theorem}[1D Dirichlet resolvent]\label{thm:Resolvent1D_upperdiag_Dir}
    Let $p\in (1,\infty)$, $k\in \NN_0$, $\tilde{\gam}\in (-1, 2p-1)\setminus\{p-1\}$, $\omega\in (0,\pi)$ and let $X$ be a $\UMD$ Banach space. Then for all $f\in W^{k,p}(\R_+, w_{\tgam+kp};X)$ and $\lambda\in \Sigma_{\pi-\omega}$ there exists a unique $u\in W^{k+2,p}_\Dir(\R_+, w_{\tilde{\gam}+kp};X)$ to $\lambda u -\delDir u =f$. Moreover, this solution satisfies
    \begin{equation}\label{eq:thm:Resolvent1D_upperdiag_Dir}
        \sum_{j=0}^2|\lambda|^{1-\frac{j}{2}}\|\d_x^j u\|_{W^{k,p}(\R_+, w_{\tgam+kp};X)}\leq C h_{k,\tgam}(\lambda) \|f\|_{W^{k,p}(\R_+, w_{\tgam+kp};X)},
    \end{equation}
    where $h_{k,\tgam}(\lambda)$ is defined in \eqref{eq:def_hlambda} and the constant $C>0$ only depends on $p,k,\tgam, \omega$ and $X$.
\end{theorem}
For the Neumann Laplacian there is a similar result, but with a shift in the smoothness of the Sobolev spaces.
\begin{theorem}[1D Neumann resolvent]\label{thm:Resolvent1D_upperdiag_Neu}
    Let $p\in (1,\infty)$, $k\in \NN_0\cup\{-1\}$, $\tgam\in (-1, 2p-1)\setminus\{p-1\}$ such that $\tgam+kp>-1$, $\omega\in (0,\pi)$ and let $X$ be a $\UMD$ Banach space. Then for all $f\in W^{k+1,p}(\R_+, w_{\tgam+kp};X)$ and $\lambda\in \Sigma_{\pi-\omega}$ there exists a unique $u\in W^{k+3,p}_\Neu(\R_+, w_{\tilde{\gam}+kp};X)$ to $\lambda u -\delNeu u =f$. Moreover, this solution satisfies
    \begin{equation}\label{eq:thm:Resolvent1D_upperdiag_Neu}
        \sum_{j=0}^2|\lambda|^{1-\frac{j}{2}}\|\d_x^j u\|_{W^{k+1,p}(\R_+, w_{\tgam+kp};X)}\leq C h_{k+ 1,\tgam}(\lambda) \|f\|_{W^{k+1,p}(\R_+, w_{\tgam+kp};X)},
    \end{equation}
    where $h_{k+1,\tgam}(\lambda)$ is defined in \eqref{eq:def_hlambda} and the constant $C>0$ only depends on $p,k,\tgam, \omega$ and $X$.
\end{theorem}

To prove Theorems \ref{thm:Resolvent1D_upperdiag_Dir} and \ref{thm:Resolvent1D_upperdiag_Neu}, we first prove some auxiliary results. We collect some estimates on the kernel $E_{\bc}$.  
\begin{lemma}\label{lem:upper-kernel}
Let $\omega\in(0,\pi)$ and $\bc\in \{\Dir, \Neu\}$. Then there exist $c,C>0$ only depending on $\omega$ such that for all $\lambda\in\Sigma_{\pi-\omega}$ and $x,y>0$ one has
\begin{equation*}
    |E_{\bc}(\lambda,x,y)|
    \leq C
    |\lambda|^{-\frac12}e^{-c|\lambda|^{\frac12}|x-y|}\cdot \begin{cases}
        \min\big\{1,|\lambda|^{\frac12 }x,|\lambda|^{\frac 12}y\big\}\quad &\mbox{if }\bc=\Dir,\\
        1\quad &\mbox{if }\bc=\Neu.
    \end{cases}
\end{equation*}
\end{lemma}
\begin{proof}
    Write $\lambda=|\lambda|e^{i\theta}$ with
$|\theta|<\pi-\omega$. For the principal square root, we have
\[
    \operatorname{Re}\sqrt\lambda
    =|\lambda|^{\frac 12}\cos(\tfrac{\theta}{2})
    \ge \sin(\tfrac{\omega}{2})|\lambda|^{\frac 12}.
\]
For $\bc=\Neu$, the result follows directly from
the definition of $E_{\Neu}$ and $x+y\ge |x-y|$.
For $\bc=\Dir$, suppose first that $x\ge y$. Then
\[
    E_{\Dir}(\lambda, x,y)
    =\frac{e^{-\sqrt\lambda(x-y)}}{2\sqrt\lambda}
      \big(1-e^{-2\sqrt\lambda y}\big).
\]
Since
\[
    |1-e^{-2\sqrt\lambda y}|
    \lesssim \min\{1,|\lambda|^{\frac 12}y\},
\]
we obtain the result. The case $y\ge x$ follows by symmetry.
\end{proof}
The following lemma is a variant of Hardy's inequality which can, for example, be derived from \cite{Mu72b}. For completeness we give an elementary proof below.
\begin{lemma}\label{lem:upper-hardy}
Let $p\in (1,\infty)$, $m\in \NN_0$, $j\in \{0,1\}$, $c>0$ and let $X$ be a Banach space. Assume that 
\begin{equation*}
    mp-1< \delta< (m+j+1)p-1,\qquad \delta \neq (j+1)p-1.
\end{equation*}
    Then for all $\eta>0$ and $f\in W^{m,p}(\R_+, w_{\delta};X)$ one has
    \begin{align*}
        \Big\|x\mapsto e^{-c\eta x}\int_0^xy^j\|f(y)\|_X \dd y\Big\|_{L^p(\R_+, w_\delta)}&\lesssim \eta^{-(j+1)}\big(1+ \eta^{-\frac{(\delta-(j+1)p+1)_+}{p}}\big)\|f\|_{W^{m,p}(\R_+, w_\delta;X)},\\
        \Big\|x\mapsto x^j\int_x^\infty e^{-c\eta y}\|f(y)\|_X\dd y\Big\|_{L^p(\R_+, w_\delta)}&\lesssim \eta^{-(j+1)}\|f\|_{W^{m,p}(\R_+, w_\delta;X)},
    \end{align*}
    where the implicit constants only depend on $p,m,j,c, \delta$ and $X$.
\end{lemma}

\begin{proof}
For any $R>0$ we find by Hölder's inequality and Hardy's inequality (see, e.g., \cite[Lemma 3.2]{LV18} using $\delta>mp-1$) that
\begin{equation}\label{eq:Hardy1}
    \begin{aligned}
    \int_0^R y^j \|f(y)\|_X\dd y &\leq \|f\|_{L^p(\R_+,w_{\delta-mp};X)}\Big(\int_0^R y^{\frac{(j+m)p-\delta}{p-1}}\dd y\Big)^{\frac{1}{p'}}\\
    &\lesssim R^{m+j+1-\frac{\delta+1}{p}}\|f\|_{W^{m,p}(\R_+, w_{\delta};X)}.
\end{aligned}
\end{equation}
Note that $m+j+1-\frac{\delta+1}{p}>0$ since $\delta< (m+j+1)p-1$. Let $R_\eta:=\min\{1, \eta^{-1}\}$. Then 
\begin{equation*}
    \int_0^x y^j \|f(y)\|_X \dd y \leq \int_0^{R_\eta} y^j \|f(y)\|_X \dd y + \one_{(R_\eta,\infty)}(x)\int_{R_\eta}^x y^j \|f(y)\|_X \dd y.
\end{equation*}
After multiplying by $e^{-c\eta x}$ and taking the $L^p(\R_+, w_\delta)$-norm it suffices to estimate the two terms on the right-hand side separately. Since the first term is independent of $x$, we find using \eqref{eq:Hardy1} that
\begin{equation*}
    \begin{aligned}
     \Big\|x\mapsto e^{-c\eta x}
       \int_0^{R_\eta}y^j\|f(y)\|_X\dd y\Big\|_{L^p(\R_+, w_\delta)} &=\Big(\int_0^\infty e^{-cp\eta x}x^\delta\dd x\Big)^\frac{1}{p} \int_0^{R_\eta} y^j \|f(y)\|_X\dd y\\
       &\lesssim \eta^{-\frac{\delta+1}{p}}R_\eta^{m+j+1-\frac{\delta+1}{p}}\|f\|_{W^{m,p}(\R_+, w_{\delta};X)}\\
     &\lesssim \min\big\{\eta^{-\frac{\delta+1}{p}}, \eta^{-(m+j+1)}\big\}\|f\|_{W^{m,p}(\R_+, w_\delta;X)}\\
     &\leq \eta^{-(j+1)}\big(1+ \eta^{-\frac{(\delta-(j+1)p+1)_+}{p}}\big)\|f\|_{W^{m,p}(\R_+, w_\delta;X)}.
\end{aligned}
\end{equation*}
For the second term, we obtain by Minkowski's integral inequality and Hölder's inequality that
\begin{equation}\label{eq:hardy2}
    \begin{aligned}
     \Big\|x\mapsto e^{-c\eta x}
       \int_{R_\eta}^x& y^j\|f(y)\|_X\dd y\Big\|_{L^p((R_\eta,\infty), w_\delta)}\\ &\leq \int_{R_\eta}^\infty y^j \|f(y)\|_X \Big(\int_y^\infty e^{-cp\eta x}x^\delta\dd x\Big)^{\frac 1p}\dd y\\
       &\leq \|f\|_{L^p((R_\eta,\infty), w_\delta;X)}\Big(\int_{R_\eta}^\infty y^{\frac{jp-\delta}{p-1}}\Big(\int_{y}^\infty e^{-cp\eta x}x^\delta \dd x \Big)^{\frac{1}{p-1}}\dd y \Big)^{\frac{1}{p'}}.
\end{aligned}
\end{equation}
By a substitution $x\mapsto \eta^{-1}x$ and $y\mapsto \eta^{-1}y$, the latter integral simplifies to
\begin{align*}
    \int_{R_\eta}^\infty y^{\frac{jp-\delta}{p-1}}\Big(\int_{y}^\infty e^{-cp\eta x}x^\delta \dd x \Big)^{\frac{1}{p-1}} \dd y&= \eta^{-(j+1)p'}\int_{\eta R_\eta}^\infty y^{\frac{jp-\delta}{p-1}}\Big(\int_{y}^\infty e^{-cp x}x^\delta \dd x \Big)^{\frac{1}{p-1}} \dd y.
\end{align*}
Since
\begin{equation*}
    \int_y^\infty e^{-cpx}x^\delta \dd x \lesssim \begin{cases}
        1\quad &\mbox{if }0<y\leq 1\\
        e^{-\frac12 cp y} &\mbox{if }y>1
    \end{cases},
\end{equation*}
and $\eta R_\eta=\min\{\eta,1\}\leq1$, we find that 
\begin{align*}
    \eta^{-(j+1)p'}\int_{\eta R_\eta}^\infty y^{\frac{jp-\delta}{p-1}}\Big(\int_{y}^\infty e^{-cp x}x^\delta \dd x \Big)^{\frac{1}{p-1}} \dd y&\lesssim \eta^{-(j+1)p'}\Big( 1+ \int_{\eta R_\eta}^1  y^{\frac{jp-\delta}{p-1}} \dd y\Big)\\
    &\lesssim \eta^{-(j+1)p'}\Big(1+ (\eta R_\eta)^{-p'\frac{(\delta-(j+1)p+1)_+}{p}}\Big)\\
    &\lesssim \eta^{-(j+1)p'}\Big(1+ \eta^{-p'\frac{(\delta-(j+1)p+1)_+}{p}}\Big).
\end{align*}
Inserting the above into \eqref{eq:hardy2} gives
\begin{align*}
    \Big\|x\mapsto e^{-c\eta x}
       \int_{R_\eta}^x& y^j\|f(y)\|_X\dd y\Big\|_{L^p((R_\eta,\infty), w_\delta)}\\ 
       & \lesssim  \eta^{-(j+1)}\big(1+ \eta^{-\frac{(\delta-(j+1)p+1)_+}{p}}\big) \|f\|_{W^{m,p}(\R_+, w_\delta;X)} 
\end{align*}
This completes the proof of the first estimate in the statement of the lemma. For the second estimate, by Minkowski's integral inequality and Hölder's inequality, we obtain
\begin{align*}
 \Big\|x\mapsto x^j
       \int_x^\infty e^{-c\eta y}\|f(y)\|_X\dd y\Big\|_{L^p(\R_+, w_\delta)}
 &\le
 \int_0^\infty e^{-c\eta y}\|f(y)\|_X
 \left(\int_0^y x^{jp+\delta}\dd x\right)^{\frac 1p}\dd y
\\
 &\lesssim
 \int_0^\infty e^{-c\eta y}\|f(y)\|_X
 y^{j+\frac{\delta+1}{p}}\dd y\\
 &\le
 \|f\|_{L^p(\R_+,w_\delta;X)}
 \Big(
 \int_0^\infty e^{-cp'\eta y}
 y^{(j+\frac{1}{p})p'}\dd y
 \Big)^{\frac{1}{p'}}\\
 &\lesssim \eta^{-(j+1)} \|f\|_{L^p(\R_+,w_\delta;X)},
\end{align*}
where for the calculation of the last integral over $y$ we have used that 
\begin{equation*}
     \tfrac{(j+1/p)p'+1}{p'}
 =j+\tfrac1p+\tfrac1{p'}=j+1.
\end{equation*}
Combining the above estimates completes the proof of the lemma.
\end{proof}

\begin{remark}\label{rem:logest} Let $m=1$ and $\delta=(j+1)p-1$.
By inspection of the proof of Lemma \ref{lem:upper-hardy}, we also have the following estimate
\begin{equation*}
 \Big\|x\mapsto e^{-c\eta x}\int_0^x
        y^j \|f(y)\|_X\,\dd y\Big\|_{L^p(\R_+,w_{\delta})}
 \lesssim
 \eta^{-(j+1)}\big[1+\log(1+\eta^{-1})\big]^{\frac{1}{p'}}
 \|f\|_{W^{1,p}(\R_+,w_{\delta};X)}.
\end{equation*}
Indeed, for $j=1$ this follows from the same proof as above and
\[
    \int_{\min\{\eta,1\}}^1\frac{\dd y}{y}
    \lesssim 1+\log(1+\eta^{-1}),
\]
in the calculation to estimate \eqref{eq:hardy2}. If $j=0$, then the estimate \eqref{eq:Hardy1} is no longer valid because Hardy's inequality fails. Instead, this estimate should be replaced by 
\begin{equation*}
    \int_0^R\|f(y)\|_X\dd y \lesssim R\big(1+\log(R^{-1})\big)^\frac{1}{p'}\|f\|_{W^{1,p}(\RR_+, w_{p-1};X)},\qquad 0< R\leq 1,
\end{equation*}
which follows from the fundamental theorem of calculus and Hölder's inequality. The rest of the proof is similar as before.
\end{remark}

The following proposition is key to establishing the precise blow-up of the resolvent in $\lambda$. 
    \begin{proposition}\label{prop:Lp_blowup}
Suppose that either of the following conditions holds.
\begin{enumerate}[(i)]
    \item Let $\bc=\Dir$, $k_\Dir:=0$ and assume that the conditions from Theorem \ref{thm:Resolvent1D_upperdiag_Dir} hold.
    \item Let $\bc=\Neu$, $k_\Neu:=1$ and assume that the conditions from Theorem \ref{thm:Resolvent1D_upperdiag_Neu} hold.
\end{enumerate}
    Then for all $\lambda\in \Sigma_{\pi-\omega}$ and $f\in W^{k+k_\bc,p}(\R_+, w_{\tgam+kp};X)$ it holds that
    \begin{equation*}
    |\lambda|\|R_{\bc}(\lambda)f\|_{L^p(\mathbb R_+,w_{\tgam+kp};X)}
    \leq C h_{k+k_\bc,\tgam}(\lambda)
    \|f\|_{W^{k+k_\bc,p}(\R_+,w_{\tgam+kp};X)}
\end{equation*}
where $h_{k+k_\bc,\tgam}(\lambda)$ is defined in \eqref{eq:def_hlambda} and the constant $C>0$ only depends on $p,k,k_\bc, \tgam, \omega$ and $X$.
\end{proposition}
\begin{proof}
By density it suffices to consider $f\in \Cc^\infty([0,\infty);X)$. We split the integral 
\begin{equation*}
    R_\bc(\lambda)f(x)  =\Big(\int_0^{x/2}+\int_{x/2}^{2x}+\int_{2x}^\infty\Big)E_{\bc}(\lambda, x,y)f(y)\dd y=:I_1+I_2+I_3,
\end{equation*}
    and consider the three regions separately. Note that by Lemma \ref{lem:upper-kernel} we have the estimate
    \begin{equation}\label{eq:upper-unified-kernel}
    |E_{\bc}(\lambda, x,y)|
    \lesssim
    |\lambda|^{-\frac12}
    e^{-c_\omega|\lambda|^{\frac12}|x-y|}
    \min\big\{1,|\lambda|^{\frac12}x,
|\lambda|^{\frac12}y\big\}^{1-k_\bc},\qquad x,y>0.
\end{equation}

    \textit{Step 1: estimate on $I_1$.} If $y\in (0, \frac{x}{2})$, then $|x-y|\geq \frac{x}{2}$ and $\min\{1, |\lambda|^\half x, |\lambda|^\half y\}^{1-k_\bc}\leq |\lambda|^{(1-k_\bc)/2}y^{1-k_\bc}$.  By \eqref{eq:upper-unified-kernel}, we find that
    \begin{equation*}
        \int_0^{\frac{x}{2}} |E_{\bc}(\lambda, x,y)|\|f(y)\|_X\dd y \lesssim |\lambda|^{-\frac{k_\bc}{2}}e^{-\frac{c}{2}|\lambda|^\half x}\int_0^xy^{1-k_\bc}\|f(y)\|_X\dd y.
    \end{equation*}
    If $\bc=\Dir$ or $\bc=\Neu$ with $\tgam\in (p-1,2p-1)$, then Lemma \ref{lem:upper-hardy} (applied to $m=k+k_\bc$, $j=1-k_\bc$, $\delta=\tgam+kp$ and $\eta=|\lambda|^\half$) gives
    \begin{equation*}
        \|I_1\|_{L^p(\R_+, w_{\tgam+kp};X)}\lesssim |\lambda|^{-1}h_{k+k_{\bc},\tgam}(\lambda) \|f\|_{W^{k+k_\bc,p}(\R_+, w_{\tgam+kp};X)}.
    \end{equation*}
    Similarly, for $\bc=\Neu$ with $\tgam\in(-1, p-1)$, we can apply Lemma \ref{lem:upper-hardy} with $m=k$ instead to obtain 
\begin{equation*}
        \|I_1\|_{L^p(\R_+, w_{\tgam+kp};X)}\lesssim  |\lambda|^{-1}h_{k+1,\tgam}(\lambda) \|f\|_{W^{k+1,p}(\R_+, w_{\tgam+kp};X)}.
    \end{equation*}
    
    \textit{Step 2: estimate on $I_2$.} Note that for $y\in (\frac{x}{2}, 2x)$ we have that $x^\theta \eqsim y^\theta$ for all $\theta\in \RR $. By \eqref{eq:upper-unified-kernel} and Young's convolution inequality, we obtain
    \begin{align*}
        \|I_2\|_{L^p(\R_+, w_{\tgam+kp};X)}&\lesssim |\lambda|^{-\half}\Big(\int_0^\infty\Big(\int_0^\infty y^{\frac{\tgam+kp}{p}}e^{-c|\lambda|^{\frac 12}|x-y|}\|f(y)\|_X\dd y\Big)^p\dd x\Big)^{\frac 1p}\\
        &\leq |\lambda|^{-\half}\Big(\int_\RR e^{-c|\lambda|^{\frac 12}|t|}\dd t \Big)\|f\|_{L^p(\R_+, w_{\tgam+kp};X)}\\
        &\lesssim |\lambda|^{-1}\|f\|_{W^{k+k_\bc,p}(\R_+, w_{\tgam+kp};X)}.
    \end{align*}

    \textit{Step 3: estimate on $I_3$.} If $y\in (2x, \infty)$, then $|x-y|\geq \frac{y}{2}$ and $\min\{1, |\lambda|^\half x, |\lambda|^\half y\}^{1-k_\bc}\leq |\lambda|^{(1-k_\bc)/2}x^{1-k_\bc}$.  By \eqref{eq:upper-unified-kernel}, we find that
    \begin{equation*}
        \int_{2x}^{\infty} |E_{\bc}(\lambda, x,y)|\|f(y)\|_X\dd y \lesssim |\lambda|^{-\frac{k_\bc}{2}} x^{1-k_\bc} \int_x^\infty e^{-\frac{c}{2}|\lambda|^\half y}\|f(y)\|_X\dd y.
    \end{equation*}
    An application of Lemma \ref{lem:upper-hardy} as in Step 1, gives the estimate
    \begin{equation*}
        \|I_3\|_{L^p(\R_+, w_{\tgam+kp};X)}\lesssim |\lambda|^{-1} \|f\|_{W^{k+k_\bc,p}(\R_+,w_{\tgam+kp};X)}.
    \end{equation*}
    Combining all the estimates for $I_1$, $I_2$ and $I_3$ yields the result.
\end{proof}

We can now prove the main theorems with sharp resolvent estimates for the Dirichlet and Neumann resolvent operators.
\begin{proof}[Proofs of Theorems \ref{thm:Resolvent1D_upperdiag_Dir} and \ref{thm:Resolvent1D_upperdiag_Neu}]
Let $k_\Dir:= 0$ and $k_\Neu: = 1$. If $\bc=\Dir$, then let $f\in \Cc^\infty(\R_+;X)$ and if $\bc=\Neu$, then let $f\in C_{\rm{c}, 1}^\infty(\overline{\R_+};X)$. Then by Lemma \ref{lem:smoothRHS} (see also \cite[Lemmas 5.3 \& 5.5]{LLRV24}) the resolvent equation $\lambda u-\del_\bc u=f $ has a unique solution given by the restricted Fourier multiplier operator
\begin{equation*}
    u_\bc (x) =  \Big(\mc{F}^{-1}\Big[\xi \mapsto \frac{ (\mc{F}\mc{E}_\bc f)(\xi)}{\lambda+ |\xi|^2}\Big]\Big)\Big|_{\R_+} =: \big(\mc{F}^{-1}[\xi \mapsto m(\xi)(\mc{F}\mc{E}_\bc f)(\xi)]\big)\big|_{\R_+},
\end{equation*}
where $\mc{E}_{\Dir}$ and $\mc{E}_\Neu$ are the odd and even extensions from $\R_+$ to $\R$, respectively. Using that $(\mc{F}^{-1}m)(x)=(2\sqrt{\lambda})^{-1}e^{-\sqrt{\lambda}|x|}$ and properties of the reflection, we find
\begin{equation*}
    u_\bc (x) = \int_0^\infty E_{\bc}(\lambda,x,y)f(y)\dd y = R_\bc(\lambda) f(x),\qquad x>0.
\end{equation*}
Therefore, by the proofs of \cite[Propositions 5.4 \& 5.6]{LLRV24} it suffices to prove \eqref{eq:thm:Resolvent1D_upperdiag_Dir} and \eqref{eq:thm:Resolvent1D_upperdiag_Neu} with $u = R_\bc(\lambda) f$ and $f$ as specified at the beginning of the proof. 
Moreover, it suffices to prove \eqref{eq:thm:Resolvent1D_upperdiag_Dir} and \eqref{eq:thm:Resolvent1D_upperdiag_Neu} with $j=0$ since the estimates for $j=2$ then follow from the resolvent equation and the estimates for $j=1$ follow by interpolation (similarly to \eqref{eq:interpolation} below). The rest of the proof is devoted to the proofs of these estimates.

From \cite[Propositions 5.4 \& 5.6]{LLRV24} we obtain for all $\lambda\in \Sigma_{\pi-\om}$ with $|\lambda|\geq 1$ that
\begin{equation*}
    |\lambda|
    \|R_{\bc}(\lambda)f\|_
      {W^{k+k_\bc,p}
       (\mathbb R_+,w_{\tgam+kp};X)}
    \lesssim
    \|f\|_{W^{k+k_\bc,p}
       (\mathbb R_+,w_{\tgam+kp};X)}.
\end{equation*}
    Therefore, it remains to prove for all $\lambda \in \Sigma_{\pi-\om}$ with $|\lambda|\leq 1$ that
    \begin{equation}\label{eq:est_smalllambda}
        |\lambda|\|R_{\bc}(\lambda)f\|_{W^{k+k_\bc,p}(\R_+, w_{\tgam+kp};X)}\lesssim h_{ k+k_\bc,\tgam}(\lambda)\|f\|_{W^{k+k_\bc,p}(\R_+, w_{\tgam+kp};X)}.
    \end{equation}
To prove \eqref{eq:est_smalllambda}, note that by Proposition \ref{prop:Lp_blowup} we obtain
\begin{align*}
    \|\d_x^2 R_\bc(\lambda)f\|_{L^p(\R_+, w_{\tgam+kp};X)}&\leq |\lambda|\|R_\bc(\lambda)f\|_{L^p(\R_+, w_{\tgam+kp};X)}+\|f\|_{L^p(\R_+, w_{\tgam+kp};X)}\\
    &\lesssim h_{ k+k_\bc,\tgam}(\lambda)\|f\|_{W^{k+k_\bc,p}(\R_+, w_{\tgam+kp};X)}.
\end{align*}
In addition, by interpolation (see \cite[Proposition 6.6]{Ro25} and a scaling argument), we have
\begin{equation}\label{eq:interpolation}
    \begin{aligned}
    |\lambda|^\half\|&\d_x R_\bc(\lambda)f\|_{L^p(\R_+, w_{\tgam+kp};X)}\\&\lesssim |\lambda|^\half \|\d_x^2 R_\bc(\lambda)f\|^\half_{L^p(\R_+, w_{\tgam+kp};X)}\|R_\bc(\lambda)f\|^\half_{L^p(\R_+, w_{\tgam+kp};X)}\\
    &\lesssim \|\d_x^2 R_\bc(\lambda)f\|_{L^p(\R_+, w_{\tgam+kp};X)}+|\lambda|\|R_\bc(\lambda)f\|_{L^p(\R_+, w_{\tgam+kp};X)}\\
    &\lesssim h_{ k+k_\bc,\tgam}(\lambda)\|f\|_{W^{k+k_\bc,p}(\R_+, w_{\tgam+kp};X)}.
\end{aligned}
\end{equation}
Since $|\lambda|\leq 1$, combining the above gives for $n\in \{0,1,2\}$ that
\begin{equation*}
    |\lambda|\|\d_x^n R_\bc(\lambda)f\|_{L^p(\R_+, w_{\tgam+kp};X)}\lesssim h_{ k+k_\bc,\tgam}(\lambda)\|f\|_{W^{k+k_\bc,p}(\R_+, w_{\tgam+kp};X)}.
\end{equation*}
For $n\in \{1, \dots, k+k_\bc-2\}$, differentiating the resolvent equation we find $\d_x^{n+2} u = \lambda \d_x^n u -\d_x^n f$ and 
\begin{equation*}
    |\lambda|\|\d_x^{n+2}R_\bc(\lambda)f\|_{L^p(\R_+, w_{\tgam+kp};X)}\leq |\lambda|^2 \|\d_x^n R_\bc(\lambda)f\|_{L^p(\R_+, w_{\tgam+kp};X)}+ |\lambda|\|\d_x^n f\|_{L^p(\R_+, w_{\tgam+kp};X)}.
\end{equation*}
Therefore, by iteration over the odd and even derivatives separately, and using that $|\lambda|\leq 1$ we find  
\begin{equation*}
    |\lambda|\|\d_x^n R_\bc(\lambda)f\|_{L^p(\R_+, w_{\tgam+kp};X)}\lesssim h_{ k+k_\bc,\tgam}(\lambda)\|f\|_{W^{k+k_\bc,p}(\R_+, w_{\tgam+kp};X)},
\end{equation*}
for all $n\in \{0, \dots, k+k_\bc\}$, which proves \eqref{eq:est_smalllambda}.
\end{proof}

\subsection{The case of small weight exponents}\label{subsec:1Dresovent_lower} In this section, we consider the one-dimensional Laplace resolvent operators on spaces $W^{k,p}(\R_+, w_\gam;X)$ with more general weight exponents $\gam>-1$. In this case one has to impose certain compatibility conditions to obtain that the operator $\lambda-\del_\bc$ is sectorial for $\lambda>0$. \\

Recall that $k_\Dir:=0$ and $k_\Neu:=1$. In addition, for $\bc\in \{\Dir, \Neu\}$ we defined
\begin{equation}\label{eq:alphabc}
     \alpha_{\bc,\gamma}:= \max\Big\{1, \frac{\gam+k_\bc p+1}{2p}\Big\},
\end{equation}
and the spaces $W^{k,p}_{\del, \bc}$ with compatibility conditions are defined in Section \ref{sec:intro}.
The main result of this section is the following resolvent estimate.
\begin{theorem}[1D resolvent with compatibility conditions]\label{thm:Resolvent1D_lowerdiag}
    Let $p\in (1,\infty)$, $k\in \NN_0$, $\bc\in \{\Dir, \Neu\}$, $\gam\in (-1, (k+2-k_\bc)p-1)\setminus\{jp-1:j\in \NN_1\}$, $\omega\in (0,\pi)$ and let $X$ be a $\UMD$ Banach space. Then for all $f\in W^{k,p}_{\del, \bc}(\R_+, w_{\gam};X)$ and $\lambda\in \Sigma_{\pi-\omega}$ there exists a unique $u\in W^{k+2,p}_{\del, \bc}(\R_+, w_{\gam};X)$ to $\lambda u -\del_\bc u =f$. Moreover, this solution satisfies
    \begin{equation}\label{eq:thm:Resolvent1D_lowerdiag}
        \sum_{j=0}^2|\lambda|^{1-\frac{j}{2}}\|\d_x^j u\|_{W^{k,p}(\R_+, w_{\gam};X)}\leq C (1+|\lambda|^{1-\alpha_{\bc, \gam}}) \|f\|_{W^{k,p}(\R_+, w_\gam;X)},
    \end{equation}
    where $\alpha_{\bc, \gam}$ is defined in \eqref{eq:alphabc} and the constant $C>0$ only depends on $p,k,\gam, \omega, \bc$ and $X$.
\end{theorem}
The core of the proof of Theorem \ref{thm:Resolvent1D_lowerdiag} is to establish the precise blow-up of $\lambda$ in the resolvent estimate. In the subsequent propositions, we will establish the precise resolvent estimates for small and large $|\lambda|$ separately. The proof of Theorem \ref{thm:Resolvent1D_lowerdiag}  is given at the end of this section.

\subsubsection{Resolvent estimates for $|\lambda|\leq 1$} We prove estimates for the resolvent $R_\bc(\lambda)$ in the case that $|\lambda|\leq 1$.

\begin{proposition}\label{prop:smalllambda}
    Let $p\in (1,\infty)$, $k\in \NN_0$, $\bc\in \{\Dir, \Neu\}$, $\gam\in (-1, (k+2-k_\bc)p-1)\setminus\{jp-1:j\in \NN_1\}$, $\om\in (0,\pi)$ and let $X$ be a $\UMD$ Banach space. Then for all $f\in W^{k,p}(\R_+, w_{\gam};X)$ and $\lambda\in \Sigma_{\pi-\om}$ with $|\lambda|\leq 1$ it holds that
    \begin{equation*}
        \sum_{j=0}^2|\lambda|^{1-\frac j2}\|\d_x^j R_\bc(\lambda)f\|_{W^{k,p}(\R_+, w_{\gam};X)}\leq C |\lambda|^{1-\alpha_{\bc,\gam}}\|f\|_{W^{k,p}(\R_+, w_{\gam};X)},
    \end{equation*}
    where $\alpha_{\bc, \gam}$ is defined in \eqref{eq:alphabc} and the constant $C>0$ only depends on $p,k,\gam, \om, \bc$ and $X$.
\end{proposition}
\begin{proof}
\textit{Step 1: large weight exponents.}
    If $\gam\in ((k-k_\bc)p-1, (k+2-k_\bc)p-1)$, then the resolvent estimates follow from Theorems \ref{thm:Resolvent1D_upperdiag_Dir} and \ref{thm:Resolvent1D_upperdiag_Neu} (applied to $k-1$ instead of $k$ for the Neumann case) since $|\lambda|\leq 1$ implies that 
    \begin{equation*}
        h_{k+k_\bc, \gam-kp}(\lambda)= 1+ |\lambda|^{1-\alpha_{\bc,\gam}}\leq 2 |\lambda|^{1-\alpha_{\bc,\gam}}.
    \end{equation*}
    \textit{Step 2: small weight exponents.} Let $\gam<(k-k_\bc)p-1$ and let $j_*$ be as in \eqref{eq:jstar} and set $\ell:=k-j_*$. Then $(\ell-k_\bc)p-1<\gam<(\ell+2-k_\bc)p-1$. Thus we can apply the result for large weight exponents from Step 1 with $\ell$ instead of $k$.  
     In particular, this gives the estimate 
     \begin{equation*}
        \sum_{j=0}^2 |\lambda|^{1-\frac{j}{2}}\|\d_x^j R_\bc(\lambda)f\|_{L^p(\R_+, w_\gam;X)}\lesssim |\lambda|^{1-\alpha_{\bc,\gam}}\|f\|_{W^{k,p}(\R_+, w_\gam;X)}.
     \end{equation*}
     Differentiating the resolvent equation and iterating separately over the odd- and even-order derivatives, as in the proofs of Theorems \ref{thm:Resolvent1D_upperdiag_Dir} and \ref{thm:Resolvent1D_upperdiag_Neu}, yields control on all derivatives of order $\leq k$. 
\end{proof}

\subsubsection{Resolvent estimates for $|\lambda|\geq 1$} We will prove estimates on the resolvent $R_\bc(\lambda)$ for large $|\lambda|$. In this case, the factors of $\lambda$ in the resolvent estimates depend on the boundary conditions. Recall that we defined
\begin{equation}\label{eq:thetabc}
    \theta_{\bc, \gam}:= k- k_\bc -\frac{\gam+1}{p}.
\end{equation}

\begin{proposition}\label{prop:largelambda}
    Let $p\in (1,\infty)$, $k\in \NN_0$, $\bc\in \{\Dir, \Neu\}$, $\gam\in (-1, (k+2-k_\bc)p-1)\setminus\{jp-1:j\in \NN_1\}$, $\om\in (0,\pi)$ and let $X$ be a $\UMD$ Banach space. Let $\lambda\in \Sigma_{\pi-\om}$ with $|\lambda|\geq 1$. Then for all $f\in W^{k,p}(\R_+, w_{\gam};X)$ it holds that
    \begin{equation}\label{eq:prop:largelambda1}
        \sum_{j=0}^2|\lambda|^{1-\frac j2}\|\d_x^j R_\bc(\lambda)f\|_{W^{k,p}(\R_+, w_{\gam};X)}\leq C |\lambda|^{\frac 12 \max\{0,\theta_{\bc,\gamma}\}}\|f\|_{W^{k,p}(\R_+, w_{\gam};X)},
    \end{equation}
    where $\theta_{\bc, \gam}$ is defined in \eqref{eq:thetabc}.
    If $f\in W^{k,p}_{\del, \bc}(\R_+, w_{\gam};X)$, then it holds that
    \begin{equation}\label{eq:prop:largelambda2}
        \sum_{j=0}^2|\lambda|^{1-\frac j2}\|\d_x^j R_\bc(\lambda)f\|_{W^{k,p}(\R_+, w_{\gam};X)}\leq C\|f\|_{W^{k,p}(\R_+, w_{\gam};X)}.
    \end{equation}
    In both cases, the constant $C>0$ only depends on $p,k,\gam, \om, \bc$ and $X$.
\end{proposition}
\begin{remark}
    The estimate \eqref{eq:prop:largelambda2} can also be obtained as a special case of Proposition \ref{prop:sect_nonoptimal}. 
\end{remark}
The estimate \eqref{eq:prop:largelambda2} will allow us to show that the (shifted) Laplace operator is sectorial on $W^{k,p}_{\del, \bc}(\R_+, w_{\gam};X)$, while it might be more natural to consider the space $W^{k,p}(\R_+, w_{\gam};X)$ without compatibility conditions. However, sectoriality on this larger space fails, since the $|\lambda|$-factor in \eqref{eq:prop:largelambda1} is
optimal. This is the content of the following proposition.
\begin{proposition}\label{prop:largelambda_optimal}
    Let $p\in (1,\infty)$, $k\in \NN_0$, $\bc\in \{\Dir, \Neu\}$, $\gam\in (-1, (k-k_\bc)p-1)\setminus\{jp-1:j\in \NN_1\}$ and let $X$ be a $\UMD$ Banach space. Then for all $\lambda>0$ large enough it holds that
    \begin{equation*}
        \lambda\|R_\bc(\lambda)\|_{\mc{L}(W^{k,p}(\R_+, w_\gam;X))}\eqsim \lambda^{\theta_{\bc, \gam}/2}.
    \end{equation*}
\end{proposition}

Before proving Propositions \ref{prop:largelambda} and \ref{prop:largelambda_optimal}, we first prove an auxiliary result which gives an explicit expression for the commutators of derivatives and the resolvent.
\begin{lemma}\label{lem:commutator}
    Let $n,\ell\in \NN_0$ with $\ell\leq n$, $m\in \{0,1\}$, $\lambda\in \Sigma_\pi$ and let $X$ be a Banach space. For all $f\in \Cc^\infty([0,\infty);X)$, we have for almost every $x>0$ that
    \begin{equation*}
        \d_x^nR_{m}(\lambda)f(x)= \d_x^{n-\ell} R_{m+\ell}(\lambda )f^{(\ell)}(x)+ (-1)^{n-m-1}e^{-\sqrt{\lambda}x}\sum_{r=0}^{\lfloor\frac{\ell-m-1}{2}\rfloor}(\sqrt{\lambda})^{n-m-2r-2}f^{(m+2r)}(0),
    \end{equation*}
    where the sum is considered to be empty if $\ell\leq m$.
\end{lemma}
\begin{proof}
    For any $g\in\Cc^\infty([0,\infty);X) $, we obtain by \eqref{eq:prop1DkernelR} and integration by parts that
    \begin{align*}
        \d_x \int_0^\infty E_m(\lambda, x, y)g(y)\dd y &=-\int_0^\infty \big[\d_yE_{m+1}(\lambda, x,y)\big] g(y)\dd y\\
        &= \int_{0}^\infty E_{m+1}(\lambda, x, y)g'(y)\dd y + E_{m+1}(\lambda, x,0)g(0).
    \end{align*}
    By iteration of the above identity and \eqref{eq:prop1DkernelR}, we obtain that
    \begin{align*}
        \d^\ell_xR_m(\lambda)f(x)&=\int_0^\infty E_{m+\ell}(\lambda, x,y)f^{(\ell)}(y)\dd y + \sum_{r=0}^{\ell-1}\big[\d_x^{\ell-1-r}E_{m+r+1}(\lambda, x,0)\big]f^{(r)}(0)\\
        &=R_{m+\ell}(\lambda)f^{(\ell)}(x) + \sum_{r=0}^{\lfloor\frac{\ell-m-1}{2}\rfloor}\big[\d_x^{\ell-1-2r-m} \lambda^{-\half} e^{-\sqrt{\lambda}x}\big]f^{(m+2r)}(0)\\
        &= R_{m+\ell}(\lambda)f^{(\ell)}(x) + (-1)^{\ell-m-1}e^{-\sqrt{\lambda}x}\sum_{r=0}^{\lfloor\frac{\ell-m-1}{2}\rfloor} (\sqrt{\lambda})^{\ell -m-2r-2}f^{(m+2r)}(0).
    \end{align*}
    The result now follows by additionally taking $n-\ell$ derivatives.
\end{proof}

We can now prove the resolvent estimates for large $\lambda$ in Proposition \ref{prop:largelambda}.
\begin{proof}[Proof of Proposition \ref{prop:largelambda}]
    \textit{Step 1: large weight exponents.} If $\gam>(k-k_\bc)p-1$, then $\theta_{\bc, \gam}<0$ and $W^{k,p}(\R_+, w_\gam;X)= W_{\del, \bc}^{k,p}(\R_+, w_\gam;X)$. The corresponding resolvent estimates are already proved in Theorem \ref{thm:Resolvent1D_upperdiag_Dir} and \ref{thm:Resolvent1D_upperdiag_Neu} since $|\lambda|\geq 1$.

    \textit{Step 2: small weight exponents.} Let $\gam\in (-1, (k-k_\bc)p-1)$ and let $j_*\in \NN_0$ be such that $(k-j_*)p-1< \gam < (k-j_*+1)p-1$, see \eqref{eq:jstar}. Let $f\in C^\infty_{{\rm c}, j_*}(\overline{\R_+};X)$ which is dense in $W^{k,p}(\R_+, w_\gam;X)$ by Lemma \ref{lem:densityLLRV}. Fix $j\in \{0,1,2\}$, $n\in \{0,\dots, k\}$ and set $\ell_n:= (n-k+j_*)_+$. In particular, we have that $0\leq \ell_n\leq j_*$ and $0\leq n-\ell_n\leq k-j_*$. Lemma \ref{lem:commutator} gives the formula
        \begin{equation}\label{eq:largelambda1}
        \begin{aligned}
            \d_x^{n+j}R_{k_\bc}(\lambda)f(x)= &\;\d_x^{n+j-\ell_n} R_{k_\bc+\ell_n}(\lambda )f^{(\ell_n)}(x)\\
            &\;+ (-1)^{n+j-k_\bc-1}e^{-\sqrt{\lambda}x}\sum_{r=0}^{\lfloor\frac{\ell_n-k_\bc-1}{2}\rfloor}(\sqrt{\lambda})^{n+j-k_\bc-2r-2}f^{(k_\bc+2r)}(0),
        \end{aligned}
    \end{equation}
    where $k_\bc=0$ and $k_\bc=1$ for Dirichlet and Neumann boundary conditions, respectively.
    We apply the resolvent estimates for large weight exponents from Step 1 with $k-j_*$ instead of $k$. Note that we apply the Dirichlet or Neumann resolvent estimate depending on whether $k_\bc+\ell_n $ is even or odd, respectively. In particular, using that $n-\ell_n\leq k-j_*$ and $\ell_n \leq j_*$, Step 1 implies that
    \begin{equation}\label{eq:largelambda2}
    \begin{aligned}
                |\lambda|^{1-\frac j2}\|\d_x^{n+j-\ell_n} R_{k_\bc+\ell_n}(\lambda) f^{(\ell_n)}\|_{L^p(\R_+, w_\gam;X)}&\lesssim \|f^{(\ell_n)}\|_{W^{k-j_*,p}(\RR_+, w_\gam;X)}\\
                &\leq \|f\|_{W^{k,p}(\R_+, w_\gam;X)}.
    \end{aligned}
    \end{equation}
    By the trace theory for weighted Sobolev spaces \cite[Theorem 4.6]{Ro25}, we obtain
    \begin{equation*}
        \|f^{(k_\bc+2r)}(0)\|_X\lesssim \|f\|_{W^{k,p}(\R_+, w_\gam;X)},\qquad k_\bc\in \{0,1\},\, r\in \{0, \dots,\lfloor\tfrac{\ell_n-k_\bc-1}{2}\rfloor \}.
    \end{equation*}
    Note that \cite[Theorem 4.6]{Ro25} applies since $k_\bc+2r \leq j_*-1< k -\frac{\gam+1}{p}$. Furthermore, since $|e^{-\sqrt{\lambda}x}|\leq e^{-c |\lambda|^\half x}$, a straightforward calculation shows that
   \begin{equation*}
       \|e^{-\sqrt{\lambda}\,\cdot}\|_{L^p(\R_+, w_\gam)}\lesssim |\lambda|^{-\frac{\gam+1}{2p}}.
   \end{equation*}
   This implies that
   \begin{equation}\label{eq:largelambda3}
   \begin{aligned}
              |\lambda|^{1-\frac j2}&|\lambda|^{\frac{n+j-k_\bc}{2}-r-1}\|e^{-\sqrt{\lambda}\, \cdot}\|_{L^p(\R_+, w_\gam)}\|f^{(k_\bc+2r)}(0)\|_X\\&\lesssim |\lambda|^{\frac{n-k_\bc}{2}-\frac{\gam+1}{2p}-r}\|f\|_{W^{k,p}(\R_+, w_\gam;X)}
              \lesssim |\lambda|^{\frac{\theta_{\bc,\gam}}{2}}\|f\|_{W^{k,p}(\R_+, w_\gam;X)},
   \end{aligned}
   \end{equation}
   where for the latter inequality we have used that $n\leq k$, $r\geq 0$ and $|\lambda|\geq 1$. Combining \eqref{eq:largelambda1}, \eqref{eq:largelambda2} and \eqref{eq:largelambda3}, and summing over $n$ and $j$ yields \eqref{eq:prop:largelambda1}.

   To prove \eqref{eq:prop:largelambda2}, we note that by Lemma \ref{lem:density} it suffices to take $f\in C_{\del, \bc, j_*}^\infty(\overline{\R_+};X)$. Therefore, all the boundary terms in \eqref{eq:largelambda1} vanish and \eqref{eq:largelambda2} implies the required resolvent estimate \eqref{eq:prop:largelambda2}.
\end{proof}

Finally, we also give the proof of Proposition \ref{prop:largelambda_optimal} showing that the factor $\lambda^{\theta_{\bc,\gam}/2}$ is sharp.
\begin{proof}[Proof of Proposition \ref{prop:largelambda_optimal}] It suffices to prove the lower estimate. Since $\gam<(k-k_\bc)p-1$, we have $\theta_{\bc,\gam}>0$ and $j_*\geq k_\bc+1$. Fix $x_0\in X$ with $\|x_0\|_X=1$, choose $\chi\in\Cc^\infty([0,\infty))$ such that $\chi=1$ on $[0,1]$, and set
\begin{equation*}
    f(x):=x^{k_\bc}\chi(x)x_0.
\end{equation*}
Thus $f^{(k_\bc)}(0)=x_0$ and $f^{(k_\bc+2r)}(0)=0$ for every $r\geq1$. Applying Lemma \ref{lem:commutator} with $n=k$ and $\ell=j_*$ gives
\begin{equation}\label{eq:optimal-commutator}
\begin{aligned}
 \d_x^kR_\bc(\lambda)f
  ={}&\d_x^{k-j_*}R_{k_\bc+j_*}(\lambda)f^{(j_*)}\\
   &+(-1)^{k-k_\bc-1}(\sqrt\lambda)^{k-k_\bc-2}
     e^{-\sqrt\lambda\,\cdot}x_0.
\end{aligned}
\end{equation}
Set $\ell:=k-j_*$, i.e., $\ell p-1<\gam<(\ell+1)p-1$.  
Proposition \ref{prop:largelambda}, applied with smoothness $\ell$ and Dirichlet or Neumann boundary condition depending on whether $k_\bc+j_*$ is even or odd, respectively, yields
\begin{equation*}
 \lambda
 \|\d_x^\ell R_{k_\bc+j_*}(\lambda)f^{(j_*)}\|_{L^p(\R_+,w_\gam;X)}
 \lesssim \|f^{(j_*)}\|_{W^{\ell,p}(\R_+,w_\gam;X)}.
\end{equation*}
On the other hand,
\begin{equation*}
 \|e^{-\sqrt\lambda\,\cdot}\|_{L^p(\R_+,w_\gam)}
 =c_{p,\gam}\lambda^{-\frac{\gam+1}{2p}}.
\end{equation*}
Therefore, \eqref{eq:optimal-commutator} and the reverse triangle inequality give
\begin{align*}
    \lambda\|R_\bc(&\lambda)f\|_{W^{k,p}(\R_+,w_\gam;X)}\\
    &\geq
    \lambda\|\d_x^kR_\bc(\lambda)f\|_{L^p(\R_+,w_\gam;X)}
    \\
    &\geq
    \lambda(\sqrt\lambda)^{k-k_\bc-2}
    \|e^{-\sqrt\lambda\,\cdot}x_0\|_{L^p(\R_+,w_\gam;X)}
    -\lambda
    \big\|\d_x^\ell R_{k_\bc+j_*}(\lambda)f^{(j_*)}
    \big\|_{L^p(\R_+,w_\gam;X)}
    \\
    &\geq
    c_{p,\gam}
    \lambda^{1+\frac{k-k_\bc-2}{2}
                 -\frac{\gam+1}{2p}}
    -C\|f^{(j_*)}\|_{W^{\ell,p}(\R_+,w_\gam;X)}
    \\
    &=c_{p,\gam}\lambda^{\theta_{\bc,\gam}/2}-C_f.
\end{align*}
This is bounded from below by $c\lambda^{\theta_{\bc,\gam}/2}$ for some $c>0$ and for all sufficiently large $\lambda$. 
\end{proof}

To conclude this section, we give the proof of Theorem \ref{thm:Resolvent1D_lowerdiag} by combining the resolvent estimates in Propositions \ref{prop:smalllambda} (small $|\lambda|$) and \ref{prop:largelambda} (large $|\lambda|$).
\begin{proof}[Proof of Theorem \ref{thm:Resolvent1D_lowerdiag}]
    It suffices to consider $\gam\in (-1, (k-k_\bc)p-1)$ since for $\gam >(k-k_\bc)p-1$ the result follows from Theorems \ref{thm:Resolvent1D_upperdiag_Dir} and \ref{thm:Resolvent1D_upperdiag_Neu}. Take $f\in C^\infty_{\del, \bc, j_*}(\overline{\RR_+};X)$, then by Lemma \ref{lem:smoothRHS} the resolvent equation $\lambda u -\del_\bc u =f $ has a unique solution $u$ which satisfies $u^{(k_\bc+2r)}(0)=0$ if $k_\bc + 2r< k+2-\frac{\gam+1}{p}$. Moreover, this solution coincides with $R_\bc(\lambda)f$. By combining the estimates in Propositions \ref{prop:smalllambda} and \ref{prop:largelambda}, we obtain \eqref{eq:thm:Resolvent1D_lowerdiag}. A density argument similarly to the proof of \cite[Proposition 5.4]{LLRV24} using Lemma \ref{lem:density} proves existence of a solution $u\in W^{k+2,p}_{\del,\bc}(\R_+, w_\gam;X)$. Finally, the uniqueness follows from Proposition \ref{prop:sect_nonoptimal}. 
\end{proof}

\section{Growth of the semigroup}\label{sec:semigroup}
In this section, we will prove Theorem \ref{thm:intro_main_thm}\ref{it:thm1D3} and \ref{it:thm1D4} about the optimal growth of the Dirichlet and Neumann heat semigroups. To this end, we will proceed in two steps.
\begin{enumerate}[(1)]
    \item\label{it:step1} Using the one-dimensional resolvent estimates from Theorem \ref{thm:Resolvent1D_lowerdiag}, we can prove Theorem \ref{thm:intro_main_thm} for $d=1$.
    \item\label{it:step2}  Using the factorisation of the semigroup in a normal and tangential part, we extend the one-dimensional result from Step \ref{it:step1} to $d\geq 2$. 
\end{enumerate}
Step \ref{it:step1} and \ref{it:step2} are carried out in Sections \ref{subsec:semigroup1D} and \ref{subsec:semigroupdD}, respectively.

\subsection{The proof of Theorem \ref{thm:intro_main_thm}\ref{it:thm1D3} and \ref{it:thm1D4} for \texorpdfstring{$d=1$}{d=1}}\label{subsec:semigroup1D}
We provide the proof of the main theorem in the special case that $d=1$. We use the resolvent estimates from Section \ref{sec:1D} and the implication \ref{it:lem:growth1}$\implies$\ref{it:lem:growth2} in Lemma \ref{lem:growth_semigroup_abstract} to get sharp estimates for the one-dimensional heat semigroups.

\begin{proof}[Proof of Theorem \ref{thm:intro_main_thm}\ref{it:thm1D3} and \ref{it:thm1D4} for $d=1$] We prove the two statements separately. 

\textit{Proof of \ref{it:thm1D3}.}
If $\gam< (2-k_\bc)p-1$, then $\alpha_{\bc,\gam}=1$. So by Theorem \ref{thm:Resolvent1D_lowerdiag}, we obtain that $-\del_\bc$ is sectorial, and boundedness of $T^1_\bc$ follows from Lemma \ref{lem:growth_semigroup_abstract}. To obtain that $-\del_\bc$ also has a bounded $\Hinf$-calculus, we argue similarly to the proof of Theorem \ref{thm:intro_main_thm}\ref{it:thm1calculus} in Section \ref{subsec:calculus}. Indeed, the
lower-order realisation $B:=-\del_\bc$ occurring in the proof of
\ref{it:thm1calculus} has a bounded $\Hinf$-calculus by \cite[Theorems 1.1 \& 1.2]{LLRV24}. Let $j_*>k_\bc$. Then repeating the similarity argument from the proof of \ref{it:thm1calculus}, with
\(
    J:=(1+B)^N
\) and $D(J)=D(B^N)$, gives that 
\(
    -\del_\bc=J^{-1}BJ
\) has a bounded
$\Hinf$-calculus of angle zero. For $j_*\leq k_\bc$ the functional calculus results for $-\del_\bc$ are already contained in \cite[Theorems 1.1 \& 1.2]{LLRV24}.

\textit{Proof of \ref{it:thm1D4}.}
If $\gam>(2-k_\bc)p-1$, then Theorem \ref{thm:Resolvent1D_lowerdiag} and Lemma \ref{lem:growth_semigroup_abstract} imply that 
\begin{equation}\label{eq:1Destimate}
    \|T^1_\bc(z)\|\lesssim 1+ |z|^{\alpha_{\bc, \gam}-1}, \qquad z\in \Sigma_{\eta},
\end{equation}
for any $\eta\in (0,\frac \pi 2)$. 
It remains to prove the lower bound. Define $\zeta\in C^{\infty}([0,\infty))$ such that $\zeta(x)=1$ for $x\leq \half$ and $\zeta(x)=0$ for $x\geq \frac 34$. Let $x_0\in X$ with $\|x_0\|_X=1$ be fixed and set $\zeta_\bc(y):=y^{1-k_\bc}\zeta(y)$. Then $\zeta_\bc^{(k_\bc+2r)}(0)=0$ for all $r\in \NN_0$ and thus $\zeta_\bc x_0\in W^{k,p}_{\del, \bc}(\R_+, w_\gam;X)$. By similar calculations as in the proof of \cite[Theorems 6.1(ii) and 6.2(ii)]{LLRV24}, we find
\begin{equation}\label{eq:lowerest1D}
    \|T^1_\bc(t)\zeta_\bc x_0\|_{L^p(\R_+, w_\gam;X)}\gtrsim t^{\frac{\gam+k_\bc p+1}{2p}-1},\qquad t\geq 1.
\end{equation}
Since the semigroup is not bounded, the semigroup property implies that $\|T^1_\bc(t)\|\geq 1$ for all $t\geq 0$.
\end{proof}
The following proposition proves the necessity of the compatibility conditions in the Sobolev spaces for sectoriality of the (shifted) Laplace operator. In the unweighted setting, this was already proved in \cite{DD11}.
\begin{proposition}\label{prop:notsect}
     Let $p\in (1,\infty)$, $k\in \NN_0$, $\bc\in \{\Dir, \Neu\}$, $\gam\in (-1, (k-k_\bc)p-1)\setminus\{jp-1:j\in \NN_1\}$ and let $X$ be a $\UMD$ Banach space. Let $\del_\bc$ be the Laplacian on $W^{k,p}(\RR_+, w_\gam;X)$ with domain
     \begin{equation*}
         D(\del_\bc):= W^{k+2,p}_\bc(\RR_+, w_\gam;X).
     \end{equation*}
     Then for every $\lambda>0$, the operator $\lambda-\del_\bc$ is not sectorial.
\end{proposition}
\begin{proof}
Fix $\lambda>0$ and assume that $A:=\lambda-\del_\bc$ is sectorial of some angle. Then, in particular, we have the estimate
\begin{equation*}
    \sup_{\mu>0} \mu\|(\mu+A)^{-1}\|_{\mc{L}(W^{k,p}(\RR_+, w_\gam;X))}<\infty.
\end{equation*}
On the other hand, $(\mu+A)^{-1} = R_\bc(\mu+\lambda)$, and Proposition \ref{prop:largelambda_optimal} gives, for all sufficiently large $\mu$, that
\begin{align*}
    \mu\|(\mu+A)^{-1}\|_{\mc{L}(W^{k,p}(\RR_+, w_\gam;X))} &= \frac{\mu}{\mu+\lambda} (\mu+\lambda)\|R_\bc(\mu+\lambda)\|_{\mc{L}(W^{k,p}(\RR_+, w_\gam;X))}\\
    &\gtrsim \frac{\mu}{\mu+\lambda} (\mu+\lambda)^{\frac{\theta_{\bc,\gam}}{2}}\to \infty \quad \text{ as }\mu\to \infty,
\end{align*}
since $\gam< (k-k_\bc)p-1$ implies $\theta_{\bc,\gam}>0$. This contradicts sectoriality of $A$.
\end{proof}

\subsection{The proof of Theorem \ref{thm:intro_main_thm}\ref{it:thm1D3} and \ref{it:thm1D4} for arbitrary dimensions}\label{subsec:semigroupdD}
We prove the main theorem for arbitrary $d\geq 1$ by using the factorisation $T_\bc^d = T^1_\bc T^{d-1}$ and the sharp semigroup estimates from above for the case $d=1$. In addition, by using the implication \ref{it:lem:growth2}$\implies$\ref{it:lem:growth1} in Lemma \ref{lem:growth_semigroup_abstract} we also obtain resolvent estimates for the Laplace operator in arbitrary dimensions.

\begin{proof}[Proof of Theorem \ref{thm:intro_main_thm}\ref{it:thm1D3} and \ref{it:thm1D4} for $d\geq 2$] We first prove the semigroup properties stated in \ref{it:thm1D3} and \ref{it:thm1D4}.
By \cite[Proposition 6.8]{Ro25} and standard properties of the trace, it holds that \begin{equation}\label{eq:part0}
    W^{k,p}_{\del, \bc}(\RRdh, w_\gam;X)= L^p(\RR^{d-1}; W^{k,p}_{\del,\bc}(\R_+, w_{\gam};X))\cap W^{k,p}(\RR^{d-1}; L^p(\R_+, w_\gam;X)),
\end{equation}
up to equivalent norms. Moreover, we note that the heat semigroup on $\RRdh$ factorises as $T_\bc^d(z) = T^1_\bc(z)T^{d-1}(z)$, where $T^{d-1}$ is the standard heat semigroup on $\RR^{d-1}$. Therefore, it suffices to estimate the two norms on the right-hand side separately.

Let $\eta\in (0, \frac \pi 2)$. Since $T^{d-1}$ is uniformly bounded on $L^p(\RR^{d-1}; Y)$ for any Banach space $Y$, and by the one-dimensional estimate \eqref{eq:1Destimate}, we obtain
\begin{equation}\label{eq:part1}
    \begin{aligned}
    \|T^d_\bc(z)f\|_{L^p(\RR^{d-1}; W^{k,p}(\R_+, w_{\gam};X))}&=\|T^{d-1}(z)T^1_\bc(z) f\|_{L^p(\RR^{d-1}; W^{k,p}(\R_+, w_{\gam};X))}\\&\lesssim \|T^1_\bc(z)f\|_{L^p(\RR^{d-1}; W^{k,p}(\R_+, w_{\gam};X))}\\&\lesssim (1+ |z|^{\alpha_{\bc,\gam}-1}) \|f\|_{L^p(\RR^{d-1}; W^{k,p}(\R_+, w_{\gam};X))},
\end{aligned}
\end{equation}
for all $z\in \Sigma_\eta$ and $f\in W^{k,p}_{\del, \bc}(\RRdh, w_\gam;X)$. Using that $\d_{\tilde{x}}^{\tilde{\alpha}}T^{d-1}(z)g= (\d^{\tilde{\alpha}}G^{d-1}_z)\ast g$ for any $\tilde{\alpha}\in \NN_0^{d-1}$, we find with Young's convolution inequality that
\begin{equation*}
    \|T^{d-1}(z)g\|_{W^{k,p}(\RR^{d-1};Y)}\lesssim \sum_{|\tilde{\alpha}|\leq k} |z|^{-|\tilde{\alpha}|/2}\|g\|_{L^p(\RR^{d-1};Y)}\lesssim \|g\|_{L^p(\RR^{d-1};Y)},
\end{equation*}
for any Banach space $Y$, $z\in \Sigma_\eta$ with $|z|\geq 1$ and $g\in L^p(\RR^{d-1};Y)$. Again using the one-dimensional estimate \eqref{eq:1Destimate}, we find
\begin{equation}\label{eq:part2}
\begin{aligned}
         \|T^d_\bc(z)f\|_{W^{k,p}(\RR^{d-1}; L^p(\R_+, w_\gam;X))}&=\|T^{d-1}(z)T^1_\bc(z) f\|_{W^{k,p}(\RR^{d-1}; L^p(\R_+, w_\gam;X))}\\
         &\lesssim \|T^1_\bc(z)f\|_{L^p(\RR^{d-1}; L^p(\R_+, w_\gam;X))}\\ &\lesssim (1+ |z|^{\alpha_{\bc,\gam}-1}) \|f\|_{L^p(\RR^{d-1}; W^{k,p}(\R_+, w_{\gam};X))},
\end{aligned}
\end{equation}
for all $z\in \Sigma_\eta$ with $|z|\geq 1$. 
Combining \eqref{eq:part1} and \eqref{eq:part2} with \eqref{eq:part0} yields 
\begin{equation*}
     \|T^d_\bc(z)f\|_{W^{k,p}(\RRdh, w_\gam;X)}\lesssim (1+ |z|^{\alpha_{\bc,\gam}-1}) \|f\|_{W^{k,p}(\RRdh, w_\gam;X)}, 
\end{equation*}
for all $z\in \Sigma_\eta$ with $|z|\geq 1$ and $f\in W^{k,p}_{\del, \bc}(\RRdh, w_\gam;X)$. Let $\sigma\in (\eta, \frac{\pi}{2})$ and since $(T_\bc^d(z))_{z\in\Sigma_\sigma}$
is an analytic $C_0$-semigroup on $W_{\del,\bc}^{k,p}(\RRdh,w_\gam;X)$, we find that 
\[
\sup_{z\in\Sigma_\eta,\, |z|\leq1}
\|T_\bc^d(z)\|<\infty.
\]
This proves the claimed upper bound for the semigroup.

We continue with the lower bound  on the semigroup in \ref{it:thm1D4}. Take $0\neq g\in \mc{S}(\RR^{d-1})$ and set $g_t(\tilde{x}):= t^{-\frac{d-1}{2p}}g(t^{-\half}\tilde{x})$ for $t>0$ and $\tilde{x}\in \RR^{d-1}$. Let $\zeta_\bc$ be as in the proof of Theorem \ref{thm:intro_main_thm} for $d=1$ and fix $x_0\in X$ with $\|x_0\|=1$. Then $h_t:=\zeta_\bc\otimes g_tx_0 \in W^{k,p}_{\del, \bc}(\RRdh, w_\gam;X)$ for all $t\geq 1$ and $\sup_{t\geq 1}\|h_t\|_{W^{k,p}(\RRdh, w_\gam;X)}<\infty$. Using the factorisation of the heat semigroup and the one-dimensional result \eqref{eq:lowerest1D}, we obtain for $t\geq 1$ that
\begin{align*}
    \|T^d_\bc(t)h_t\|_{L^p(\RRdh, w_\gam;X)} &= \|T^1_\bc(t)\zeta_\bc\|_{L^p(\RR_+, w_\gam)}\|T^{d-1}(t) g_tx_0\|_{L^p(\RR^{d-1};X)}\\
    &=\|T^1_\bc(t)\zeta_\bc\|_{L^p(\RR_+, w_\gam)}\|T^{d-1}(1) g x_0\|_{L^p(\RR^{d-1};X)}
    \gtrsim t^{\frac{\gam+k_\bc p+1}{2p}-1}.
\end{align*}
Hence, the semigroup $T_\bc^d$ on $W^{k,p}_{\del, \bc}(\RRdh, w_\gam;X)$ satisfies
\begin{equation*}
    \|T_\bc^d(t)\|\geq \frac{\|T^d_\bc(t)h_t\|_{W^{k,p}(\RRdh, w_\gam;X)}}{\|h_t\|_{W^{k,p}(\RRdh, w_\gam;X)}} \gtrsim  \|T^d_\bc(t)h_t\|_{L^p(\RRdh, w_\gam;X)}\gtrsim t^{\frac{\gam+k_\bc p+1}{2p}-1}.
\end{equation*}
Combining the upper and lower estimates gives the required properties for the semigroup in \ref{it:thm1D3} and \ref{it:thm1D4}. 

It remains to show that Theorem \ref{thm:intro_main_thm}\ref{it:thm1D1} and \ref{it:thm1calculus} also hold with $\lambda=0$ if $\gam<(2-k_\bc)p-1$. 
The above semigroup estimates together with Lemma \ref{lem:growth_semigroup_abstract} imply  sectoriality of $ - \del_\bc$. To see that $-\del_\bc$ also has a bounded $\Hinf$-calculus, one can argue similarly to the one-dimensional case, see the proof of Theorem \ref{thm:intro_main_thm}\ref{it:thm1D3} for $d=1$ in Section \ref{subsec:semigroup1D}.
\end{proof}

We conclude with two remarks concerning some extensions of our methods. First, for the critical weight exponent $\gam=(2-k_\bc)p-1$, we also expect to obtain a precise blow-up rate for the heat semigroups as discussed in the following remark.

\begin{remark}\label{rem:critical-log}
Let $k\geq1$ and
\(
    \gam=(2-k_\bc)p-1.
\)
Although we do not identify the domain of the generator for this weight
exponent, we conjecture that the heat semigroup satisfies
\begin{equation}\label{eq:loga}
        \|T_\bc^d(t)\|_{\mc L(
      W_{\del,\bc}^{k,p}(\RRdh,w_\gam;X))}
    \eqsim
    \bigl(1+\log(1+t)\bigr)^{\frac{1}{p'}},
    \qquad t\geq0.
\end{equation}
Indeed, using the estimate in Remark \ref{rem:logest} instead of Lemma \ref{lem:upper-hardy} in the proof of the one-dimensional resolvent estimate (Propositions \ref{prop:Lp_blowup} and \ref{prop:smalllambda}) gives
\[
    |\lambda|\|R_\bc(\lambda)\|
    \lesssim
    \bigl(1+\log(1+|\lambda|^{-1})\bigr)^{\frac{1}{p'}},
    \qquad 0<|\lambda|\leq1.
\]
Redoing the proof of
Lemma \ref{lem:growth_semigroup_abstract} (see \cite[Lemma 2.3]{LLRV24}) with a logarithmic blow-up, yields the upper bound for the semigroup, and
the extension to $d\geq2$ follows from the factorisation used in
Section~\ref{subsec:semigroupdD}.

For the lower bound, choose
$\chi_0,\chi_\infty\in C^\infty(\R_+)$ such that
\[
 \chi_0=0\text{ on }[0,1],\quad \chi_0=1\text{ on }[2,\infty),
 \qquad
 \chi_\infty=1\text{ on }[0,\tfrac12],\quad
 \chi_\infty=0\text{ on }[1,\infty),
\]
set $\chi_R(y):=\chi_0(y)\chi_\infty(y/R)$. Let $x_0\in X$ with $\|x_0\|_X=1$ and define
\[
    f_R(y):=(1+\log R)^{-\frac 1p}
    y^{-(2-k_\bc)}\chi_R(y)x_0.
\]
Then $\sup_{R>4}\|f_R\|_{W^{k,p}(\RR_+, w_\gam;X)}<\infty$, and one can show that
\[
    \|T_\bc^1(R^2)f_R\|_{L^p(\RR_+, w_\gam;X)}
    \gtrsim(1+\log R)^{\frac{1}{p'}}.
\]
In this way one can prove the precise logarithmic blow-up in \eqref{eq:loga}.\\
\end{remark}

We finish with a remark about a weighted setting where the heat semigroups are actually bounded. 
\begin{remark}\label{rem:rho}
    The growth of the heat semigroup is caused by the growth of the weight $w_\gam$ at infinity. Instead consider the following weights: let $\rho\in C^\infty([0,\infty))$ be  such that $\rho(x_1)>0$ if $x_1>0$ and 
\begin{equation*}
    \rho(x_1)=|x_1| \quad \text{ for }0\leq x_1\leq 1\quad \text{ and}\quad  \rho(x_1)=2 \quad \text{for }x_1\geq 2.
\end{equation*}
For $\gam>-1$ and $x=(x_1,\widetilde x)\in\RRdh$, define the weight    $\rho_\gam(x):=\rho(x_1)^\gam.$
Thus $\rho_\gam=w_\gam$ in a neighbourhood of $\partial\RRdh$, whereas
$\rho_\gam$ is constant for $x_1\geq2$. We define the spaces
\begin{equation*}
    L^p(\RRdh, \rho_\gam;X),\quad W^{k,p}(\RRdh,\rho_\gamma;X)\quad \text{ and }\quad  W_{\del,\bc}^{k,p}(\RRdh,\rho_\gamma;X)
\end{equation*}
similarly to Section \ref{sec:intro}. As the weights $\rho_\gam$ have the same behaviour as $w_\gam$ close to the boundary, the traces in the above Sobolev spaces are again well defined. 

Under the assumptions of Theorem~\ref{thm:intro_main_thm}, for every
$\sigma\in(0,\frac{\pi}{2})$ we expect that
\begin{equation}\label{eq:rho}
    \sup_{z\in\Sigma_\sigma}
    \|T_\bc^d(z)\|_{\mc L(
       W_{\del,\bc}^{k,p}(\RRdh,\rho_\gam;X))}
    <\infty,
\end{equation}
which can be proved with the same strategy as we have applied in Sections \ref{sec:1D} and \ref{sec:semigroup}. The part of the one-dimensional
resolvent proof which uses the behaviour of $w_\gam$ at infinity is
Lemma \ref{lem:upper-hardy}. Splitting the integrals at $x=1$, using
the original weighted Hardy estimate on $(0,1)$ and the unweighted
estimate on $(1,\infty)$, where $\rho_\gam\eqsim1$, gives
\[
\left\|
 x\mapsto e^{-c\eta x}\int_0^x y^j\|f(y)\|_X\,\dd y
\right\|_{L^p(\R_+,\rho_\delta)}
\lesssim
\eta^{-(j+1)}
\|f\|_{W^{m,p}(\R_+,\rho_\delta;X)}
\]
under the assumptions of Lemma \ref{lem:upper-hardy}. Thus the factor $|\lambda|^{1-\alpha_{\bc,\gam}}$ in Proposition~\ref{prop:smalllambda} disappears. The rest of the arguments in Sections \ref{subsec:semigroup1D} and \ref{subsec:semigroupdD} can be used with minor modifications to prove \eqref{eq:rho}.
\end{remark}

\subsection*{AI disclosure statement}
ChatGPT 5.6 by OpenAI was used to explore proof strategies for this paper. The paper was entirely written by the author, who takes full responsibility for the content.

\bibliographystyle{plain}
\bibliography{references}
\end{document}